\documentclass[acmsmall]{acmart}

\AtBeginDocument{%
  \providecommand\BibTeX{{%
    \normalfont B\kern-0.5em{\scshape i\kern-0.25em b}\kern-0.8em\TeX}}}

\setcopyright{iw3c2w3}
\copyrightyear{2020}
\acmYear{2021}
\acmDOI{10.1145/3447958}
\acmJournal{TOMS}
\acmVolume{47}
\acmNumber{3}
\acmArticle{24}
\acmMonth{9}

   \usepackage{graphicx}

    \usepackage{url}

    \usepackage{amsmath,amsfonts,amsthm}
    \usepackage[utf8]{inputenc}
    \usepackage[T1]{fontenc}
    \usepackage{ifpdf}
    \usepackage[english]{babel}
    \usepackage[labelfont=normalsize,textfont=normalsize]{subfig} 
    \usepackage{hyperref}
    \usepackage{upgreek}

    \usepackage{algorithm2e}
    
    \SetCommentSty{mycommfont}

    \usepackage{tikz}
    \usepackage{todonotes}
    
    \usepackage{xcolor,mdframed}
    
    \usetikzlibrary{arrows}

    \usepackage{tikz}
\usepackage{subfig}
\DeclareMathOperator{\real}{\text{Re}}
\DeclareMathOperator{\imag}{\text{Im}}

\newcommand{\me}{\mathrm{e}}
\newcommand{\mi}{\mathrm{i}}

\newcommand{\Mat}[1]{\ensuremath{\mathbf{#1}}}
\newcommand{\vect}[1]{\ensuremath{\mathbf{#1}}}
\newcommand{\Dict}{\Mat{D}}
\newcommand{\Gram}{\Mat{G}}
\newcommand{\atom}{\vect{d}}
\newcommand{\win}{\vect{g}}
\newcommand{\coef}{\vect{c}}
\newcommand{\sig}{\vect{x}}
\newcommand{\res}{\vect{r}}
\newcommand{\normany}[2]{ \lVert #1 \rVert_{#2} }
\newcommand{\normtwo}[1]{\normany{#1}{2}}
\newcommand{\abs}[1]{| #1 |}
\newcommand{\inprod}[2]{ \langle #1,#2 \rangle }

\usepackage{xcolor}

\definecolor{darkviolet}{rgb}{0.58,0,0.83} 

\usepackage{pgfplots}

\usepackage{datatool}
\DTLsetseparator{ = }
\DeclareMathOperator*{\argmax}{argmax} 

    \theoremstyle{plain}
    \newtheorem{theorem}{Theorem}
    \newtheorem{corollary}[theorem]{Corollary}

 \newcommand{\cres}{\ensuremath{\boldsymbol{\varrho}}}
 \newcommand{\sumcres}{\ensuremath{\rho}}
 
 \newenvironment{fundef}[1][]{%
   \begin{mdframed}[%
      backgroundcolor={gray!25}, hidealllines=true,
      skipabove=0.7\baselineskip, skipbelow=0.7\baselineskip,
      splitbottomskip=2pt, splittopskip=2pt, #1]%
   \makebox[0pt]{
      \smash{
         \fontsize{24pt}{24pt}\selectfont
         \hspace*{-19pt}
         \raisebox{-2pt}{
            {\color{gray!70!black}\sffamily\bfseries :}
         }%
      }%
   }%
}{\end{mdframed}}

\begin{document}
%
\title{Fast Matching Pursuit with Multi-Gabor Dictionaries}

\author{Zden\v{e}k Pr\r{u}\v{s}a}
\email{zprusa@kfs.oeaw.ac.at}
\author{Nicki Holighaus}
\email{nicki.holighaus@oeaw.ac.at}
\author{Peter Balazs}
\email{peter.balazs@oeaw.ac.at}
\affiliation{%
  \institution{Acoustics Research Institute, Austrian Academy of Sciences}
  \streetaddress{Wohllebengasse 12--14}
  \city{Vienna}
  \state{Austria}
  \postcode{1040}
}
\thanks{A repository with code reproducing Fig.\,\ref{fig:kernels} and the timing 
scenario in Section\,\ref{sec:mptk} is available at
\url{https://github.com/ltfat/fastmpwithmultigabor}.}

\renewcommand{\shortauthors}{Pr\r{u}\v{s}a, et al.}

\begin{abstract}
    Finding the best $K$-sparse approximation of a signal in a
redundant dictionary is an NP-hard problem.
Suboptimal greedy matching pursuit (MP) algorithms are generally
used for this task. 
In this work, we present an acceleration technique and an implementation
of the matching pursuit algorithm acting on a multi-Gabor dictionary, i.e., 
a concatenation of
several Gabor-type time-frequency dictionaries, each of which consisting of translations 
and modulations of a possibly different window and time and frequency
shift parameters. 
The technique 
is based on pre-computing and thresholding inner products between
atoms and on updating the residual directly in the coefficient domain, i.e., 
without the round-trip to the signal domain.
Since the proposed acceleration technique involves an approximate update step,
we provide theoretical and experimental results illustrating the convergence of
the resulting algorithm.
The implementation is written in C (compatible with C99 and C++11)
and we also provide Matlab and GNU Octave interfaces.
For some settings, the implementation is up to 70 times faster than the
standard Matching Pursuit Toolkit (MPTK).

\end{abstract}



\maketitle


\section{Introduction}
\label{sec:intro}
The best $K$-sparse approximation of a signal
$\sig\in\mathbb{R}^L$ in an
overcomplete dictionary of $P$ normalized atoms (vectors)
$\Dict=\left[\atom_0|\atom_1|\dots|\atom_{P-1}\right]\in\mathbb{C}^{L\times
P}$,
$\normtwo{ \atom_p } =1$
is an NP-hard problem \cite{damaav97}.
Given the budget of $K$ nonzero elements of the coefficient vector
$\coef\in\mathbb{C}^P$,
the problem can be formally written as the minimization of the approximation
error in the energy norm $\normtwo{ \sig - \Dict\coef}$ such that
\begin{equation}
    \min \normtwo{ \sig - \Dict\coef} \ \
    \text{subject to} \ \ \normany{ \coef}{0} \leq K,
    \label{eq:ksparseproblem}
\end{equation}
where the zero ``norm'' $\normany{ \cdot }{0}$ returns the number
of non-zero elements.
A similar problem is the minimization of 
$\normany{ \coef }{0}$ given the approximation error tolerance $E$
\begin{equation}
    \min \normany{ \coef}{0} \ \ \text{subject to} \ \ \normtwo{ \sig - \Dict\coef}\leq E.
\end{equation}
Both problems can be tackled by employing greedy \emph{matching pursuit} 
(MP) algorithms. The only difference is the choice of the stopping criterion.
However, greedy algorithms are known to be suboptimal in the sense that they 
are not guaranteed to choose the best combination of $K$ atoms.
Instead, an approximation rate i.e. the decrease of the approximation error with
iterations has been studied.
It has been shown that the basic version of MP \cite{mazh93}
achieves an exponential approximation rate
\cite{damaav97,dete96,grfiva06,grva06}.
To date, several variants of generic MP and its orthogonal version OMP \cite{parekr93,damazh94}
were proposed e.g.
complementary MP \cite{ragu08,ragu10}, 
cyclic MP \cite{stchr10,stchrgr11},
gradient pursuit \cite{blda08a,blda08b}, 
local OMP \cite{magrbiva09,magrvabi11} and
self projected MP \cite{rerose17}.
In practice, without imposing any structure on
the dictionary,
the effectiveness of the algorithms quickly
deteriorates when increasing the dimensionality
of the problem; either by increasing the input signal length $L$
or the size of the dictionary $P$.
Even with structured dictionaries,
which allow usage of fast algorithms in place of matrix operations,
a naive implementation can still be prohibitively inefficient;
e.g. processing even just a few seconds of an audio signal, 
which typically consist of tens of thousands 
of samples per second, can take hours.

An overview of greedy algorithms, a class of algorithms MP falls under, can be found in \cite{yada09,rish15}
and in the context of audio and music processing in \cite{st09,plbldagrda10,za16}.
Notable applications of MP algorithms include audio analysis 
\cite{grderobama96},
\cite{gr01},
coding \cite{stgi06,rarida08,chaneba14}, 
time scaling/pitch shifting 
\cite{de07}
\cite{stdacu06}, 
source separation \cite{gr02}, denoising \cite{bhde14},
partial and harmonic detection and tracking \cite{leda06}
and EEG analysis \cite{du07}.


In this contribution, we present a method for accelerating MP-based
algorithms acting on a single overcomplete Gabor dictionary or on a 
concatenation of several Gabor dictionaries with possibly different windows 
and parameters (hence the term multi-Gabor dictionary).
The main idea of the present acceleration technique is 
performing the residual update in the coefficient domain while
exploiting the locality
of the inner products between the atoms in the dictionaries
and dismissing values below a user definable threshold.
It is then feasible to store all significant inner products
in a lookup table and avoid atom synthesis and the residual re-analysis in every
iteration of MP as it is usually done in practice.
The size of the lookup table as well as the cost of computing it
are independent of the signal length and
depend only on the parameters of the Gabor dictionaries.
An integral part of this contribution is the freely available implementation
in C (compatible with C99 and C++11), which can be found in the
backend library of the Matlab/GNU Octave Large Time-Frequency Analysis Toolbox (LTFAT, \url{http://ltfat.github.io})
\cite{ltfatnote015,ltfatnote030}
available individually at \url{http://ltfat.github.io/libltfat}.
The low level C language (or rather a subset of C99 and C++11 standards)
was chosen for two reasons: 
First, the nature of the MP algorithm does not lend itself to 
an efficient implementation in a high level language due to its overhead.
For example, a proof-of-concept Matlab implementation was about 50 times slower than the
final C implementation.
Second, a C-based shared (dynamic) library can be interfaced
from most of the high and even low level languages. 
The programming interface (documentation available at
\url{http://ltfat.github.io/libltfat/group__multidgtrealmp.html}) was designed 
with this use case in mind.
Indeed, since version 2.3.0, LTFAT  itself interfaces the library trough a MEX function, whose call
is wrapped in a function \texttt{multidgtrealmp} (see
\url{http://ltfat.github.io/doc/gabor/multidgtrealmp.html}).

To date, considering a vast body of literature dealing with MP,
surprisingly few authors address effective (non-textbook) implementation of the algorithm
let alone provide code.
In the original paper,
Mallat and Zhang \cite[Appendix E]{mazh93} proposed to perform the residual update
in the coefficient domain using inner products between the atoms.
They present an analytic formula for evaluating the inner products between atoms
of a multi-scale Gabor dictionary with a Gaussian window.
%
%
An implementation by Ferrando~et.~al.~\cite{fe02} is tailored to the Gaussian window-based 
multiscale Gabor dictionary defined on an interval. The authors choose to trade 
updating the residual directly in the coefficient domain
for the flexibility in choosing the dictionary parameters and in boundary handling.
The de-facto standard implementation of several MP based algorithms is in the 
Matching Pursuit Toolkit (MPTK) \cite{krgr06}.
The toolbox is not restricted to Gabor dictionaries,
and, therefore, the coefficient-domain update rule is not exploited. 
In comparison, 
the present method and implementation is applicable to general
multi-Gabor dictionaries while being much faster than MPTK.

The paper is organized as follows. Section \ref{sec:prelim} summarizes the 
necessary theoretical background of the MP algorithm, and introduces the approximate residual update crucial to the proposed implementation, as well as a complementary convergence result.  Section \ref{sec:resupd} introduces the main contribution of the paper:
the method for accelerating MP iterations with multi-Gabor dictionaries.
The rest of the paper consists of 
Section \ref{sec:details} which discusses some practical aspects of the proposed 
method and Section \ref{sec:mptk} containing timing and approximation quality comparisons with the reference implementation in MPTK.

\section{Preliminaries} \label{sec:prelim}
Matrices will be denoted with bold capital upright letters,
e.g., $\Mat{M}$, column vectors with lowercase bold upright letter
such as $\vect{x}$. Conjugate transpose will be denoted with a star superscript,
($\vect{x}^*,\Mat{M}^*$),
scalar variables with a capital or lowercase italics letter $s,S$ and
scalar constants as upright capital or lowercase letters like $\uppi,\me,\mi$.
A single element of a matrix or a vector will be selected using round brackets
$\Mat{M}(m,n)$, $\vect{x}(l)$.
The index is always assumed to be applied modulo vector length 
(or matrix size in the respective direction)
such that  $\vect{x}(l)=\vect{x}(l + kL)$ for $l=0,\dots,L-1$ and $k\in\mathbb{Z}$.
Moreover, we will use two indices and subscript for vectors such that 
$\coef(m,n)_{M} = \coef(m + nM)$
in order to transparently ``matrixify'' a vector. 
Sub-vectors and sub-matrices will be selected by an index set denoted by a 
caligraphic letter e.g. $\vect{x}({\mathcal{P}})$ 
and the $m$-th row of a matrix $\Mat{M}$
 will be selected using the notation $\Mat{M}(m,\bullet)$ and the $n$-th column by $\Mat{M}(\bullet,n)$, respectively.
We will omit brackets when indexing the outcome of the matrix-vector or matrix-matrix product
i.e. we will use $\Mat{M}\vect{x}(p)$ instead of 
$\left(\Mat{M}\vect{x}\right)(p)$.
Scalar-domain functions 
used on matrices or vectors are applied element-wise e.g.
$\abs{\vect{x}}^2(l)=\abs{\vect{x}(l)}^2$. 
The inner product of two vectors in $\mathbb{C}^L$ is given as
$\inprod{\vect{x}}{\vect{y}} =
\vect{y}^*\vect{x} = 
\sum_{l=0}^{L-1} \vect{x}(l)\overline{\vect{y}(l)}$,
where the overline denotes complex conjugation.
Real and imaginary parts of a complex number will be denoted as $\real(c)$ 
and $\imag(c)$ respectively
and the phase as $\arg(c)$ 
such that $c=\real(c) + \mi\imag(c) = \abs{c}\me^{\mi \arg(c)}$.
The $2$--norm of a vector is defined as 
$\normany{ \vect{x}}{2} = (\sum_{l=0}^{L-1} \abs{\vect{x}(l)}^2 )^{1/2}$.
In particular, the $2$--norm relates to the inner product as $\normany{\vect{x}}{2}^2 = \inprod{\vect{x}}{\vect{x}}$. 
For a matrix $\Mat{M}$, $\|\Mat{M}\| = \|\Mat{M}\|_2 = \max_{\|\vect{x}\|_2=1} \|\Mat{M}\vect{x}\|_2$ is the matrix norm induced by the 2--norm.

%

\subsection{Multi-Gabor Dictionaries}
A Gabor dictionary $\Dict_{\left(\win,a,M\right)}$ 
generated from a window $\win\in\mathbb{R}^L,\normtwo{\win}=1$,
time shift $a$ and a number of modulations $M$
is given as 
\begin{equation}
    \begin{split}
        \Dict_{\left(\win,a,M\right)}(l,m+nM)
        & =  \win(l-na)\me^{\mi 2\uppi m (l - na) /M} 
    \label{eq:gabdict}
    \end{split}
\end{equation}
for $l=0,\dots,L-1$ and $m=0,\dots,M-1$ for each $n=0,\dots,N-1$, where 
$N=L/a$ is the number of window time shifts and
the overall number of atoms is $P=MN$.
The expression $(l-na)$ is assumed to be evaluated modulo $L$ according to
the circular indexing. 
The redundancy of a dictionary will be defined as $P/L=M/a$.
A~multi-Gabor dictionary consisting of $W$ Gabor dictionaries is defined as
\begin{equation}
    \left[ %
        \Dict_{\left(\win_1,a_1,M_1\right)}\middle|  
        \Dict_{\left(\win_2,a_2,M_2\right)}\middle|  
        \dots
        \Dict_{\left(\win_{W},a_{W},M_{W}\right)}  
        \right]
\end{equation}
and we will also use a shortened notation
$\Dict_w=\Dict_{\left(\win_w,a_w,M_w\right)}$.
Generally, $a_w,M_w$ need only be divisors of $L$. 
Due to technical reasons explained in Sec.~\ref{sec:multi}, however,
efficiency of the presented algorithm depends 
on the pairwise 
\emph{compatibility} of $a_u,a_v$ and $M_u,M_v$, implying  
some restrictions of the dictionary parameters. 
In the following, we focus on the optimal setting, i.e., 
parameters $a_w$ 
chosen such that every pair $a_u$, $a_v$ is divisible by
$a_\text{min}=\min\left\{a_u,a_v\right\}$ and, similarly,
every pair of $M_u$, $M_v$ should divide  
$M_\text{max}=\max\left\{M_u,M_v\right\}$
and each $M_w/a_w$ should be a positive integer.
While not strictly necessary, such setting
is commonly used in practice and
leads to the most efficient implementation.

\subsection{Matching Pursuit -- MP}
In this section we recall the idea behind the MP algorithm, summarize its steps
and explain an alternative way of performing the MP iterations exploiting
the inner products between the atoms.

Recall that the main goal is to find the best $k$-term approximation of a given signal $\sig$ by elements from the dictionary, i.e. $\sig \approx \sig_k = \sum \limits_{l=1}^k c_l \atom_l.$ 
The MP algorithm iteratively decreases the approximation error 
(energy of the residual) $E_{k+1}=\normtwo{\res_{k+1}}^2$ by considering orthogonal projections 
$\inprod{\res_k}{\atom_p}\atom_p$ of the residual $\res_{k} = \sig - \sig_{k}$ over the individual 
$P$ elements of the normalized dictionary. The
element, $p_\text{max}$, which decreases the energy of the residual 
$\normtwo{ \res_{k} - \inprod{\res_k}{\atom_p}\atom_p}^2$
the most is selected and the residual is updated: 
$\res_{k+1} = \res_{k} - \inprod{\res_k}{\atom_{p_\text{max}}}\atom_{p_\text{max}}$.
Since the energy of the new potential residual can be written as
\begin{equation}
    \normtwo{ \res_{k} - \inprod{\res_k}{\atom_p}\atom_p}^2 =
    \normtwo{\res_k}^2 - \abs{ \inprod{\res_k}{\atom_p} }^2,
\end{equation}
the best atom to choose is the one with the highest inner product with the residual
i.e.
\[
p_{\text{max}} = \argmax_p  \abs{ \inprod{\res_{k}}{\atom_p}}.
\]
The procedure is repeated until the desired approximation error is achieved or alternatively
some other stopping criterion is met e.g. a sparsity or a selected inner product magnitude
limits are reached.
The error is usually normalized and converted to decibels by $10\log_{10}
E_{k+1}/\normtwo{\sig}^2$. 
It is known that the matching pursuit (MP) algorithm and its derivatives
can benefit from
pre-computing inner products between the atoms in the dictionary
$\Gram(k,j)=\inprod{ \atom_j}{ \atom_k}$
i.e. from pre-computing the Gram matrix
$\Gram=\Dict^*\Dict\in\mathbb{C}^{P\times P}$. 
With $\cres_k = \Dict^*\res_{k}$ denoting the coefficient-domain residual, the residual update step can be written 
as (\cite[Ch. 12]{wtour})
\begin{equation}
    \cres_{k+1} = \cres_{k} - \coef(p_{\text{max}})\Gram(\bullet,p_{\text{max}}).
    \label{eq:update}
\end{equation}
Formally, the coefficient-domain matching pursuit algorithm is summarized in Alg.~\ref{alg:CDmp}. The stopping criterion may contain several conditions, and the algorithm terminates if any of these conditions  is met. Typical stopping conditions include reaching a certain error, a maximum number of atom selections or the largest entry in $\cres_k$ falling below some value.

\begin{algorithm}[h]
\DontPrintSemicolon
\caption{Coefficient-Domain Matching Pursuit}
\label{alg:CDmp}
\KwIn{Input signal $\sig$, dictionary Gram matrix $\Gram = \Dict^* \Dict$}
\KwOut{Solution vector $\coef$}
\textbf{Initialization:}  
$\coef =\vect{0}$, 
$\cres_0 = \Dict^* \sig$,
$E_0 = \normtwo{\sig}^2$, 
$k = 0$
\;
\While{ Stopping criterion not met }
{
    \begin{enumerate}
        \item \textbf{Selection:} \label{item:selection}
             $p_{\text{max}} \leftarrow \argmax\limits_p  \abs{ \cres_k(p) }$
        \item \textbf{Update:} \label{item:update}
        \begin{enumerate}
            \item \textbf{Solution:}  \label{item:solupdate}
                $\coef(p_{\text{max}}) \leftarrow \coef(p_{\text{max}}) + \cres_k(p_{\text{max}})$ \; 
            \item \textbf{Error:} \label{item:errupdate}
                $E_{k+1} \leftarrow E_k - \abs{ \cres_k(p_{\text{max}})}^2$\;
            \item \textbf{Residual:} \label{item:resupdate}
                $\cres_{k+1} \leftarrow \cres_{k} - \cres_k(p_{\text{max}})\Gram(\bullet,p_{\text{max}})$\;
        \end{enumerate}
    \end{enumerate}
$k \leftarrow k + 1$\;
}
\end{algorithm}
%
This modification has the advantage of removing the necessity 
of synthesizing the residual and recomputing the inner product in the selection step. On the other hand,
such approach is usually dismissed as impractical in the literature due to
the high memory requirements for storing the Gram matrix.
This is however not the case for a well behaved multi-Gabor dictionary, for which the 
Gram matrix can be precomputed and significant values stored compactly, see Section \ref{sec:resupd}.

\section{Approximate update by a truncated Gram matrix} \label{ssec:approxupdate}
As already mentioned, our acceleration technique works with the coefficient
domain update formula~\eqref{eq:update}. Whenever the elements of the dictionary $\Dict$ 
are localized, most of the entries of the Gram matrix are close to zero. Discarding these 
entries is the first step towards reducing the memory requirements of coefficient-domain matching pursuit, 
at the cost of introducing a small approximation error. Before considering the additional structure 
imposed by a multi-Gabor dictionary, we discuss the implications of using a truncated Gram matrix in 
Step \ref{item:resupdate} in Alg. \ref{alg:CDmp} and provide a basic worst-case error estimate.

For this purpose, denote by $\Gram_{\epsilon}$, $\epsilon>0$, the hard-thresholded Gram matrix, i.e., 
\begin{equation}\label{eq:truncgram}
  \Gram_{\epsilon}(k,j) = \begin{cases}
                            \Gram(k,j) & \text{ if } \abs{\Gram(k,j)} > \epsilon,\\
                            0 & \text{ otherwise.}
                          \end{cases}
\end{equation}
Similar to $\cres_k$, we further denote by 
$p_k$ the index selected in the $k$-th selection step. 
Using the full Gram matrix $\Gram$ in Alg. \ref{alg:CDmp}, we always have 
\begin{equation}\label{eq:residualequivalence}
    \cres_k = \Dict^\ast\res_k,\quad \text{i.e.,}\quad \cres_k(p) = \langle \res_k,\atom_p\rangle, \text{ for all } p,
\end{equation}
where $\res_k = \sig - \Dict \coef_k$ is the true residual, i.e., the difference between the signal $\sig$ and the proposed solution after the $k$-th step $\Dict \coef_k$. If we use the truncated Gram matrix $\Gram_\epsilon$ instead, both the error $E_{k}$ and the coefficient-domain residual $\cres_k$ is Alg. \ref{alg:CDmp} are merely estimates of the true quantities. Thus, \eqref{eq:residualequivalence} does not hold anymore. Further, the selection of $p_{\text{max}}$ is based directly on the coefficient-domain residual (estimate) $\cres_k$, such that the sequence $(p_k)_k$ of selected positions must be expected to differ between matching pursuit and this approximate variant. 

Unless noted otherwise, we will use the notation $\vect{c}_0 = \vect{0}$, $\cres_0 = \Dict^\ast \sig$ and, for $k\geq 0$,
\begin{equation}\label{eq:approxResDef}
  p_{k+1} 
  = \argmax\limits_p  \abs{ \cres_k(p)},\quad 
  \cres_{k+1} 
  = \cres_{k} - \cres_k(p_{k+1})\Gram_{\epsilon}(\bullet,p_{k+1}),\quad 
  \vect{c}_{k+1} 
  = \vect{c}_{k} + \cres_k(p_{k+1})\vect{e}_{p_{k+1}},
\end{equation}
where $\vect{e}_{j}\in\mathbb{C}^P$ is the $j$-th standard unit vector.

In the following, we show that despite these differences to accurate matching pursuit, the proposed approximate coefficient-domain matching pursuit algorithm reduces the approximation error $\|\res_k\|_2^2$, unless $\max_p |\cres_k(p)|$ is too small. 
Although arbitrarily small approximation error cannot be guaranteed when this approximate scheme is used directly, we show that the approximate matching pursuit can be nested inside a simple reset procedure to ensure $\|\res_{k}\|_2^2 \rightarrow 0$, see Alg. \ref{alg:appCDMPwReinit}. Similar to Alg. \ref{alg:CDmp}, both the reset and stopping criterion may consist of any number of conditions, terminating the loop if any condition is met. For this nested execution of approximate matching pursuit, we further provide a decay estimate on $\|\res_k\|_2^2$. The proofs of the following results can be found in the Appendix.

\begin{algorithm}[h]
\DontPrintSemicolon
\caption{Approximate Coefficient-Domain Matching Pursuit with Reset}
\label{alg:appCDMPwReinit}
\KwIn{Input signal $\sig$, truncated dictionary Gram matrix $\Gram_\epsilon$, dictionary $\Dict$}
\KwOut{Solution vector $\coef^{\text{out}}$}
\textbf{Initialization:}  
$\coef^{\text{out}} =\vect{0}$, 
$\res^{\text{out}} = \sig$,
$E^{\text{out}} = \normtwo{\sig}^2$, 
$l = 0$, $k=0$, $k^{\text{out}}=0$
\;
\While{ Stopping criterion not met }
{
    \begin{enumerate}
        \item \textbf{Approximate Matching Pursuit:} \label{item:ACDMP}\\
             \KwIn{$\res^{\text{out}}$, $\Gram_\epsilon$}
             \KwOut{Solution vector $\coef$}
             \textbf{Initialization:} $\coef = 0$, 
                                      $\cres_{k} = \Dict^* \res^{\text{out}}$,\\
                                      $E_k = E^{\text{out}}$\\
             \While{ Reset criterion not met }
             {
             \begin{enumerate}
                    \item \textbf{Selection:} \label{item:selectionA}
                        $p_{\text{max}} \leftarrow \argmax\limits_p  \abs{ \cres_k(p) }$
                    \item \textbf{Update:} \label{item:updateA}
                    \begin{enumerate}
                        \item $\coef(p_{\text{max}}) \leftarrow \coef(p_{\text{max}}) + \cres_k(p_{\text{max}})$ \; 
                        \item\label{item:errupdateA}
                        $E_{k+1} \leftarrow E_k - \abs{ \cres_k(p_{\text{max}})}^2$\;
                        \item $\cres_{k+1} \leftarrow \cres_{k} - \cres_k(p_{\text{max}}) \Gram_\epsilon(\bullet,p_{\text{max}})$\;
                    \end{enumerate}                   
            \end{enumerate}
            $k \leftarrow k + 1$\;
            }
        \item \textbf{Update:} \label{item:updateB}
        \begin{enumerate}    
            \item \textbf{Selections:} \label{item:selupdate2}
                $k^{\text{out}} \leftarrow k$\;
            \item \textbf{Solution:} \label{item:solupdate2}
                 $\coef^{\text{out}} \leftarrow \coef^{\text{out}} + \coef$\;        
            \item \textbf{Residual:} \label{item:resupdate2}
                $\res^{\text{out}} \leftarrow \res^{\text{out}} - \Dict \coef$\;
            \item \textbf{Error:} \label{item:errupdate2}
                $E^{\text{out}} \leftarrow \normtwo{\res^{\text{out}}}^2$\;                
        \end{enumerate}
    \end{enumerate}
$l \leftarrow l + 1$\;
}
\end{algorithm}

\begin{theorem}\label{thm:Resdecrease}
  Fix some positive Gramian threshold $\epsilon > 0$ 
  and some $0<\delta<1/2$. Let $p_k$ and $\cres_k$, $k\geq 0$, be as in \eqref{eq:approxResDef}. If 
  \begin{equation}\label{eq:crit1a}
    |\cres_k(p_{k+1})| \geq \frac{\epsilon}{\delta} \sum_{l=1}^{k} |\cres_{l-1}(p_{l})|\quad \text{and}\quad
    |\cres_k(p_{k+1})| >  \frac{2\delta}{1-2\delta} \|\res_k\|_2,
  \end{equation}
  then $\|\res_{k+1}\|_2^2 < \|\res_{k}\|_2^2$.
\end{theorem}
  
%

The conditions \eqref{eq:crit1a} are mostly of theoretical interest for two reasons: Firstly, the estimates made in proving Theorem \ref{thm:Resdecrease} are highly pessimistic in the sense that they assume the worst-case error. The actual error is highly likely to be significantly smaller. Secondly, verifying the second condition in \eqref{eq:crit1a} is expensive, as it requires the computation of the true time-domain residual. Therefore, while possible, we will not use \eqref{eq:crit1a} as conditions for the stopping criterion. 

In the following theorem, $\cres_k$ is used as in Algorithm \ref{alg:appCDMPwReinit}, i.e., when the reset criterion is trigered in the $(k+1)$-th selection step, then the inner loop is restarted with $\cres_k := \Dict^\ast (\sig - \Dict\coef_k)$. Otherwise, the defintion of $\cres_k$ coincides with the one given above.


\begin{theorem}\label{thm:Resdecrease2}
   Assume that $\Dict$ is finite and spans $\mathbb{C}^L$ with $\lambda_{\text{min}} := 
  \inf_{\substack{\sig\in \mathbb{C}^L\\ \|\sig\|_2 = 1}} \max_p |\langle \sig, \atom_p\rangle|^2 > 0$. Fix some $0<\delta<1/2$ and $0<\varepsilon < 1-2\delta$, such that $\frac{2\delta}{1-(2\delta+\varepsilon)} < \sqrt{\lambda_{\text{min}}}$. If the Alg. \ref{alg:appCDMPwReinit} is initialized with $\Gram_\epsilon$, $0<\epsilon<1$, and the \emph{reset criterion} includes the conditions 
  \begin{equation}\label{eq:crit1}
    |\cres_k(p_{k+1})| < \frac{\epsilon}{\delta} \sum_{l=k^{\text{out}}+1}^{k} |\cres_{l-1}(p_{l})|\quad \text{and}\quad 
    |\cres_k(p_{k+1})| < \frac{2\delta (1-\varepsilon (1+\delta)^{-2}  \lambda_{\text{min}})^{\frac{k^{\text{out}}-k}{2}}}{1-(2\delta+\varepsilon)} \normtwo{\res^{\text{out}}},
  \end{equation}
  then the residual decreases exponentially: 
  \begin{equation}
    \|\res_{k+1}\|_2^2 < (1-\varepsilon (1+\delta)^{-2} \lambda_{\text{min}}) \|\res_k\|_2^2, \text{ for all } k\geq 0.
  \end{equation}
\end{theorem}

\begin{corollary}\label{cor:convergence}
  Fix some positive $\epsilon > 0$ and $\delta,\varepsilon$ as in Theorem \ref{thm:Resdecrease2}. 
  If Algorithm \ref{alg:appCDMPwReinit} is initialized $\Gram_\epsilon$, the \emph{stopping criterion} includes the condition $E^{\textnormal{out}} \leq E \|x\|_2^2$, for some $0<E<1$, and the  \emph{reset criterion} includes the conditions \eqref{eq:crit1} 
  and
  \begin{equation}\label{eq:crit3}
   k > \frac{\log(E)}{\log(1-\varepsilon (1+\delta)^{-2} \lambda_{\text{min}})},
  \end{equation}
  then the algorithm terminates after a finite number of total selection steps, achieving the desired approximation error $E^{\textnormal{out}} \leq  E \|x\|_2^2$.
\end{corollary}

Once more, the conditions given in Theorem \ref{thm:Resdecrease2} and Corollary \ref{cor:convergence} are very conservative worst-case conditions and it is not advisable to use them in practice. This is easy to see from the proofs presented in the Appendix. Nonetheless, the results serve as justification for the introduction of resets to ensure convergence. Although the results may suggest that resets are  required rather regularly, it is in practice rarely required. In the case of the multi-Gabor dictionaries considered in this paper, excellent approximation quality is achieved without resets for $\epsilon = 10^{-4}$. Heuristically, we observed that it is sufficient to reset only after a large number of selection steps, except when the truncation threshold $\epsilon$ is large, see Section \ref{sec:practconv}, where we propose an efficient stopping condition that has proven sufficient in all our experiments.

\section{Faster Approximate Coefficient-Domain Residual Update} \label{sec:resupd}
In this section, we will show that, for multi-Gabor dictionaries, the significant values 
of the Gram matrix $\Gram$
can be precomputed at a cost independent of the entire signal length $L$, 
truncated and stored efficiently for a single as well as for a multi-Gabor dictionary.
We exploit the fact that the Gram matrix $\Gram$ of a single Gabor dictionary
$\Dict=\Dict_w$ 
is highly structured. 
Using this structure, we obtain a feasible, and in fact highly efficient, implementation of 
approximate coefficient-domain matching pursuit. 
In fact, 
$\Gram$ takes the form of a \emph{twisted convolution} \cite{chabible} matrix
with a fixed kernel
$\vect{h} = \Gram(\bullet,0) = \Dict^* \win  \in\mathbb{C}^{{M N}}$; a coefficient vector
consisting of inner products of the window with all its possible time and frequency shifts.
The $p$-th column of the Gram matrix $\Gram(\bullet,p)$ i.e., the inner products of the atom at the time-frequency position $(k,j)$ (such that $p=k+jM$) with atoms at all time-frequency positions $(m,n)$ are constructed by shifting and modulating the kernel such that
\begin{equation}\label{eq:GramKern}
    \Gram(m+nM,k+jM) = \mathbf{h}(m-k,n-j)_M
    \me^{\mi 2 \uppi k \frac{a}{M}(n-j)  }
\end{equation}
for $m=0,\dots,M-1$ for each $n=0,\dots,N-1$. Crucially, $\vect{h}(\bullet,\bullet)_M$ is essentially supported around the origin for localized, low-pass windows, such that it can be truncated and stored efficiently. By considering the concentration of the kernel relative to the time-frequency index $(k,j)$, i.e. shifting $\vect{h}(\bullet,\bullet)_M$ from the origin to $(k,j)$, we see that the modulation factor in \eqref{eq:GramKern} is in fact independent of the time index $j$. We can define
\begin{equation}\label{eq:hkj} 
\vect{h}^{(k)}(m,n) := \vect{h}(m,n)_M e^{i2\pi k \frac{a}{M}  n}.
\end{equation}
Therefore, after selecting atom $p_\text{max}=k_\text{max}+j_\text{max}M$
in the MP algorithm, the coefficient-domain residual update in the style of \eqref{eq:update} reduces to a subtraction of
a truncated, modulated and weighted kernel
from the neighborhood of the time-frequency position
$(k_\text{max},j_\text{max})$. With $\vect{h}^{(k)}$ as in \eqref{eq:hkj}, the residual update can be written as
\begin{equation}
    \label{eq:fastkernupdate}
    \cres_{k+1}(\mathcal{K},\mathcal{J})_M = \cres_{k}(\mathcal{K},\mathcal{J})_M - \cres(p_\text{max})\vect{h}^{(k_\text{max})}
\end{equation}
assuming the $\vect{h}^{(k)}$ have already been truncated and $\mathcal{K}$ and $\mathcal{J}$ denote
index sets of the appropriate neighborhood encompassing the overlay with the truncated
kernel. The neighborhood $(\mathcal{K},\mathcal{J})_M$ is of fixed size and shape, independent of $p_\text{max}$, and centered on $p_\text{max}$.

\subsection{Pre-computing the Kernel for a Single Dictionary}
\begin{figure}[!tb]
    \centering
    \includegraphics[height=4.2cm,trim={0cm 0 2.3cm 0},clip]{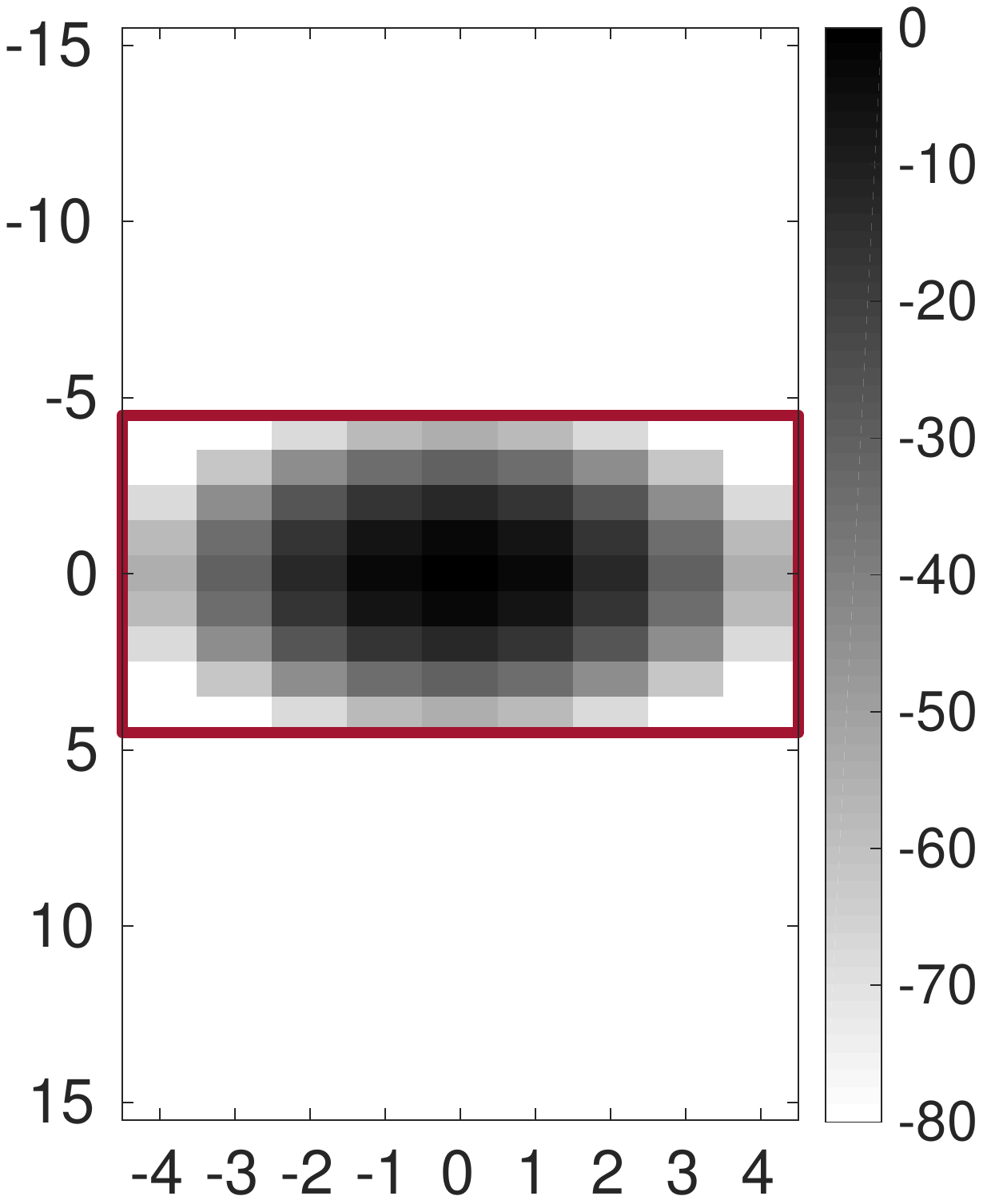}\ \
    \includegraphics[height=4.2cm,trim={0 0 2.3cm 0},clip]{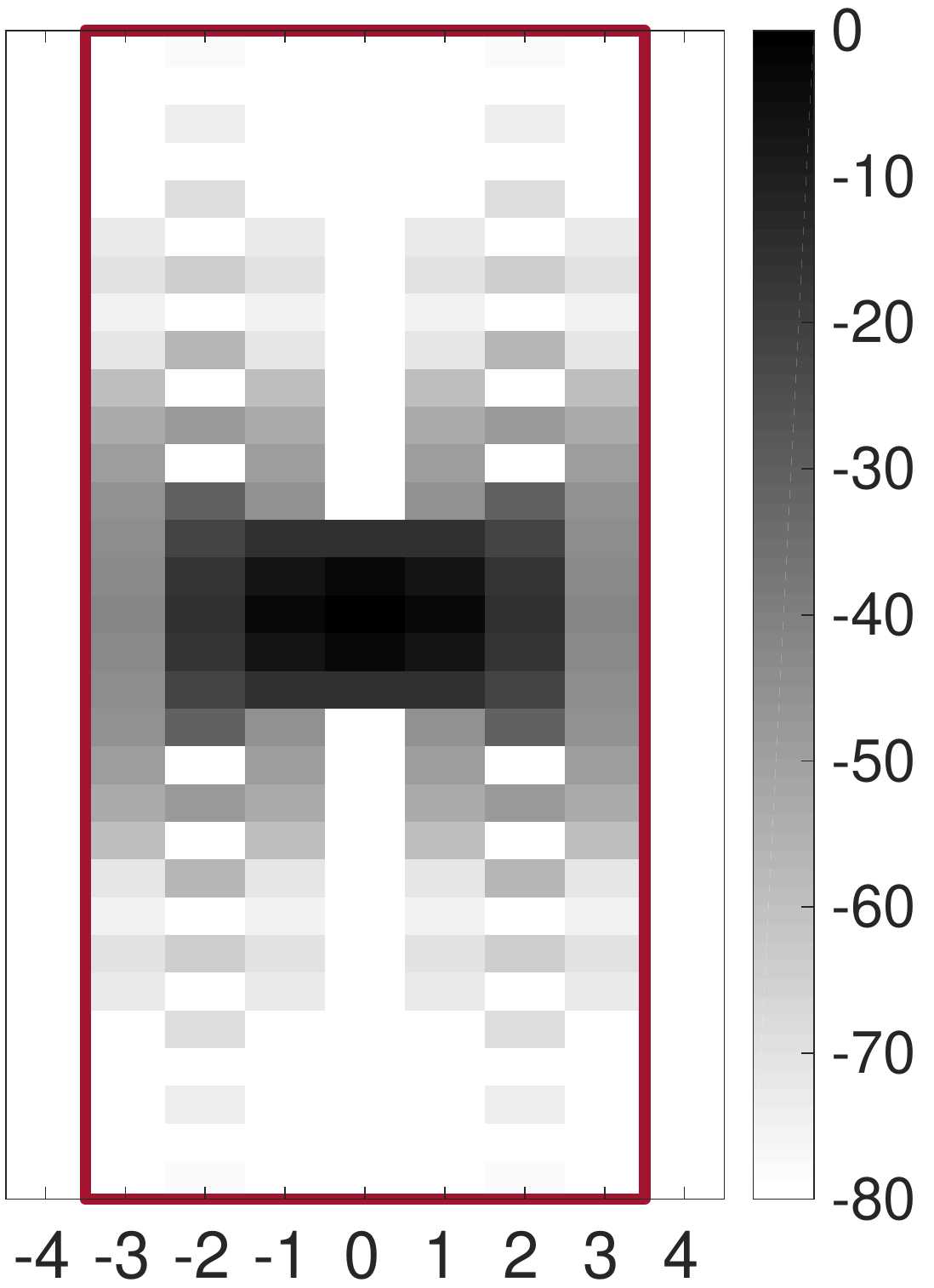}\ \
    \includegraphics[height=4.2cm]{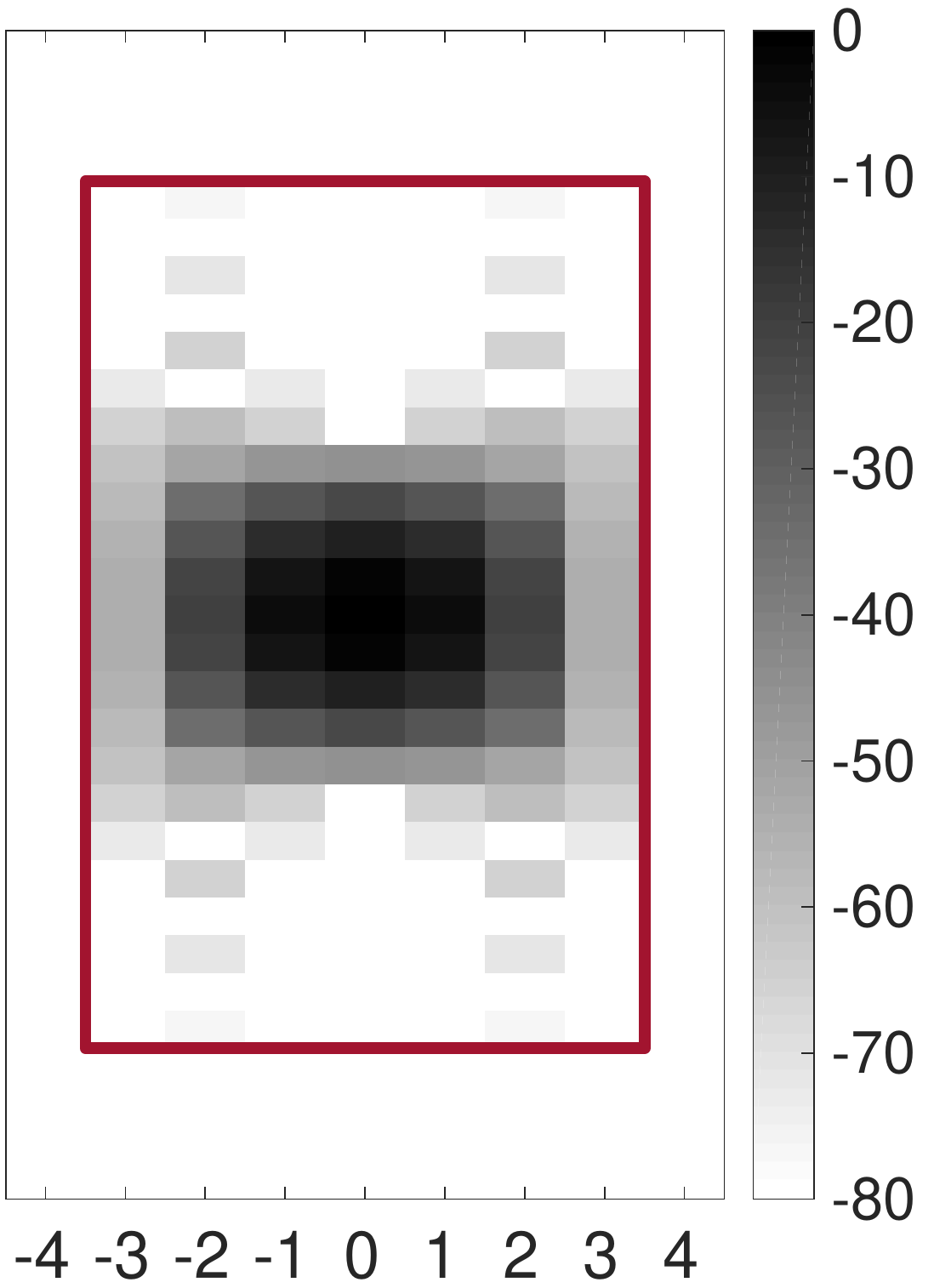}
    \caption{Examples of abs. values of truncated kernels for
    (left) Gaussian $(9\times9)$, (middle) Hann $(31\times 7)$ and (right) Blackman $(23\times 7)$  windows. }
    \label{fig:kernels}
\end{figure}

The inner products between the window and its 
translations which are sufficiently far away
are obviously zero or at least negligible.
Therefore after determining the length of the
window's effective support $L_w$, one can compute the minimum 
admissible $L_\text{min}$ for the computation of 
values of the kernel $\vect{h}$ as being the next integer multiple
of $a$ bigger than twice the length of the effective support, i.e., $L_\text{min} = \lfloor 2L_w/a \rfloor a + 1$. 
The kernel is further truncated also in the frequency direction
such that values below a certain threshold are dismissed.
The resulting size of the truncated kernel depends on the shape of 
the window and on the length of the time and frequency steps. 
The kernel size directly determines the number of
complex multiplications and additions required to perform the entire residual update step
and, obviously, also the memory requirements to store the kernel.
Examples of abs. values of kernels for several windows using time shift 
$a=512$ and $M=2048$ frequency bins are depicted in Fig.~\ref{fig:kernels}.
The values are in dB relative to the maximum element with 0~dB. Values below $-80$~dB
were cropped and are not used in the residual update. 
The threshold selection is a trade-off between faster updates (higher threshold) or less requirement for resets (lower threshold).
The idea of truncating the kernel $\vect{h}$ originates from Le~Roux~et~al.~\cite{leroux10}
who used it for replacing the operation of the (cross) Gram matrix in an iterative
scheme dealing with magnitude-only reconstruction.
The authors of the aforementioned paper also noticed that the kernel is conjugate symmetric 
about both the horizontal (time) and the vertical (frequency) axes,
which could be exploited for reducing the number of multiplications further.
%
When inspecting formula \eqref{eq:hkj}, 
it is obvious that for a fixed frequency position $k$
the modulation by $2\uppi ka/M$ radians is performed on
all rows of the kernel independently. Moreover, the modulation frequencies are
$\text{lcm}(a,M)/a$
periodic in $k$ and, therefore, all unique complex exponentials
can be tabulated and stored. In the best case when $M$ is integer divisible by $a$,
the memory requirements are equal to
storing $M/a$ additional rows of the kernel.
The cost of applying the modulation during the residual update step is
one complex multiplication per kernel column.

\subsection{Pre-computing Cross-Kernels Between Dictionaries}
\label{sec:multi}
\begin{figure}[htpb]
    \begin{center}
    \input{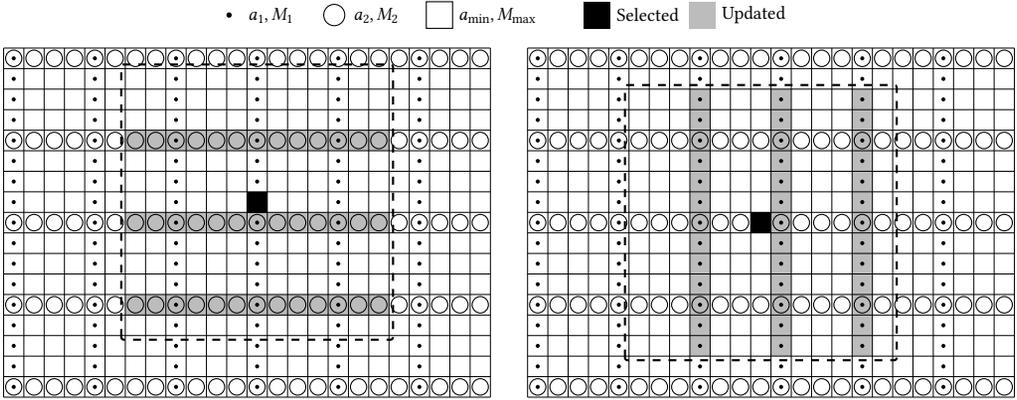}
    \end{center}
    \caption{Illustration of the residual update between dictionaries using cross kernels. The left figure shows the case when a coefficient from the dictionary 1 was selected
in the selection step of MP and
inner products with dictionary 2 are being updated and vice versa right.
The dashed line square depicts the area covered by a cross-kernel with respect to the
common grid $a_{\text{min}}, M_\text{max}$.}
    \label{fig:update}
\end{figure}

The Gram matrix of a multi-Gabor dictionary consists of Gram matrices
of individual dictionaries $\Dict_w$ and cross-Gram matrices \cite{xxlfinfram1} 
between the dictionaries.
Denoting a cross-Gram matrix as $\Gram_{w,v} = \Dict^*_{w} \Dict_{v}$
the overall Gram matrix is a block matrix with the following structure
\begin{equation}
\begin{bmatrix}
    \Gram_{1,1} &  \Gram_{1,2} & \dots  & \Gram_{1,W}\\
    \Gram_{2,1} &  \Gram_{2,2} & \dots  & \Gram_{2,W}\\
    \vdots      &  \vdots      & \ddots & \vdots\\
    \Gram_{W,1} &  \Gram_{W,2} & \dots  & \Gram_{W,W}\\
\end{bmatrix}.
\end{equation}
A cross-Gram matrix $\Gram_{w,v}$ shares the same twisted convolution 
structure with the regular Gram matrix with kernel $\vect{h}_{w,v}=\Dict^*_w \win_v$ 
only if the time-frequency shifts
are equal i.e. $a_{w}=a_{v}$ and $M_{w}=M_{v}$ .
In the case the parameters differ, the direct twisted convolution structure is lost.
The structure can be recovered on a finer ``common'' time-frequency grid given by
the time step $\text{gcd}(a_{w}, a_{v})$ and the number of frequency bins
$\text{lcm}(M_{w},M_{v})$.
The most efficient case is achieved when $a_{w}$ and $a_{v}$ are divisible by
$a_{\text{min}} = \min\left\{ a_{w}, a_{v} \right\}$
and $M_{w}$ and $M_{v}$ both divide $M_\text{max} = \max\left\{ M_{w}, M_{v} \right\}$
resulting to a common grid given by $a_{\text{min}}$ and $M_\text{max}$.
%
In the residual update step of the inner products of the residual with
the $w$-th dictionary, the modulated kernel is subsampled
by ratios $a_w/a_\text{min}$ and $M_\text{max}/M_w$ in horizontal and vertical
directions respectively.
%
To illustrate, consider a multi-Gabor dictionary consisting of two Gabor dictionaries
$\Dict_1=\Dict_{(\win_1,a_1,M_1)}$ and $\Dict_2=\Dict_{(\win_2,a_2,M_2)}$
with $a_1=4a_\text{min}, M_1=8 a_\text{min} $ and $a_2=a_\text{min}, M_2=2a_\text{min}$.
Both cross-kernels $\mathbf{h}_{1,2}$ and $\mathbf{h}_{2,1}$ 
are computed
with $a_\text{min}=a_2$ and $M_\text{max}= M_1$.
The example in Fig.~\ref{fig:update} depicts an update of inner products 
of the residual with both dictionaries on the common grid.

\section{Practical Considerations}
\label{sec:details}
Since MP is a simple algorithm, the main focus will be on a detailed description of the
accelerated residual update step from Sec.~\ref{sec:resupd}.
Other steps will be mentioned in less detail. 

The presented description is adapted to the setting of real signals $\sig\in\mathbb{R}^L$. Therefore, only the first $M_{w,\mathcal{R}}=\lfloor M_w/2 \rfloor +1$ frequency
bins form each Gabor dictionary $w$ will be considered and the 
reduced dictionary will be denoted as $\Dict_{w,\mathcal{R}}$.
A real signal can be recovered from reduced coefficient vectors $\coef_{w,\mathcal{R}}$
as
\begin{equation}
    \label{eq:idgtreal}
    \sig = \sum_{w=1}^{W}\Dict_{w,\mathcal{R}} \coef_{w,\mathcal{R}} + 
    \overline{ \Dict_{w,\mathcal{R}} } \overline{\coef_{w,\mathcal{R}}}
\end{equation}
where elements of $\overline{ \coef_{w,\mathcal{R}}}$ are set to zero whenever the conjugated partner 
is missing i.e. for frequency index $m=0$ and $m=M_w/2$ if $M_w$ is even.
Obviously, in practice the matrix operations are replaced by an efficient FFT-based
algorithm (see e.g. \cite{po76,ltfatnote011}) or by memory efficient atom-by-atom synthesis as
it is done in MPTK.

As discussed in \cite[Appendix B]{gr01}, dealing with real signals in this way requires that conjugated pairs of 
complex atoms $\atom,\overline{\atom}$ are considered as real atoms, such that the signal approximation
and the residual are real at any stage of the algorithm and the positive-negative
frequency conjugate symmetry of the coefficients is preserved. Consequently, all inner products must be adjusted as
\begin{equation}
\label{eq:dualatprod}
    \widetilde{c}=
    \frac{\inprod{\res_k}{\atom} -
    \inprod{\atom}{\overline{\atom}}
    \overline{ \inprod{\res_k}{{\atom}}}
    }{1 - \abs{\inprod{\atom}{\overline{\atom}}}^2}
\end{equation}
and subtracting the pair of atoms from the residual decreases its energy by
\begin{equation}\label{eq:projenfaster}
    \widetilde{E}=
         \left[ 
    \real(\widetilde{c})^2 +
    \imag(\widetilde{c})^2
    + \real(\inprod{\atom}{\overline{\atom}})\left(
    \real( \widetilde{c} )^2 -
    \imag( \widetilde{c} )^2
    \right) 
    - 2 
    \imag( \inprod{\atom}{\overline{\atom}} )
     \real( \widetilde{c} ) 
    \imag( \widetilde{c} )
    \right] \cdot 2.
\end{equation}
Obviously, as long as $\inprod{\atom}{\overline{\atom}}\approx 0$,
we can simplify the equations such that it is enough to consider only a single atom $\atom$.
In the following, we address the issue of an efficient search for the maximum inner product.

As discussed above, theory suggests that \eqref{eq:projenfaster} should be used to determine atoms for selection.
In practice, however, discarding the effect of the nonzero inner product between the conjugated
atoms in the selection step (setting $\inprod{\atom}{\overline{\atom}} = 0$) does not 
have a significant impact on the sparsity of the result achieving a specified approximation error
and leads to about 20\% overall speedup.
Our implementation supports both options and the technically correct one will be 
referred to as \emph{pedantic}.
On the other hand, it is crucial that the selected coefficient is adjusted using
\eqref{eq:dualatprod} prior to the residual update step.
Finally, due to the conjugate symmetry, 
the residual update step now involves a pair of atoms each of which can however
be treated separately. 
In the description of the implementation, we will work exclusively with the
reduced dictionary, therefore we will drop the ${\mathcal{R}}$ subscript in order to
lighten the notation.

Given $W$ sets of Gabor dictionary parameters $(\win_w,a_w,M_w)$ for $w=1,\dots,W$,
the initialization involves pre-computing (or loading) cross-kernels and the 
complex exponentials according to the rules described in Sec.~\ref{sec:resupd}.
The algorithm itself then starts by computing inner products of the input
signal $\res_0=\sig$ (the initial residual) with atoms from all reduced dictionaries 
$\coef_{w}
= \Dict_{w}^* \sig$ for all $w=1,\dots,W$. Whenever the algorithm is 
reset, this procedure is repeated, with the current residual $\res^{\textrm{out}}$ substituting for
 $\sig$.

\subsection{Keeping Track Of The Maximum}
\label{sec:max}
\begin{figure*}[t]
    \centering
    \includegraphics[width=14cm]{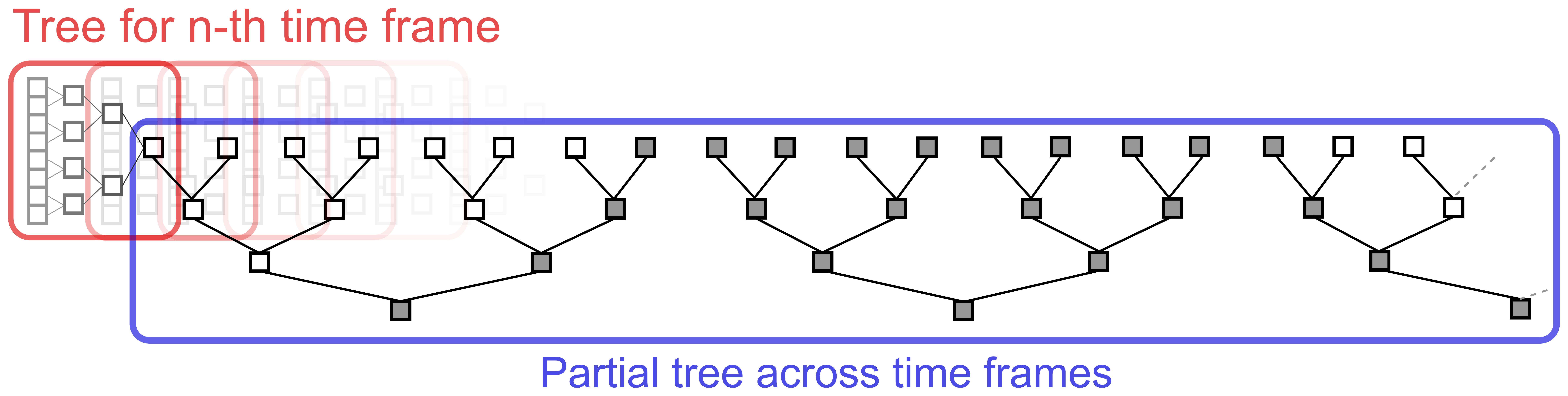}
   \caption{Schematic of the tournament tree structure for keeping track of maxima. A red border indicates individual tournament trees across frequency bins for each time frame. The blue border indicates the partial tournament tree across time frames. An example of the worst-case 10 element bottom-up 
   update of the tree across frames is shown in gray.}
    \label{fig:ttree}
\end{figure*}

Performing the full search for the maximum inner product in each selection
step is highly inefficient. 
The authors of MPTK \cite{krgr06} proposed to store positions
of maxima for each window time position and organize them
in a partial hierarchical tournament-style tree. 
Such tree contains at each level maxima from pairs from one level below.
Since the residual update affects only a limited number of 
neighboring time positions, a bottom-up tree update becomes more efficient than 
a full search.
To quantify the reduction of the number of comparisons required to find the maximum,
consider an array of length $L$ and a tree of depth $d$, where $d=0$ means a fallback to a full array search.
After the tree has been initialized, 
the search can be performed at the top level of the tree which requires only $\lceil
L/2^d\rceil-1$ comparisons.
When $Q$ consecutive elements from the array are modified, the bottom up tree update requires 
additional $Q+d$ comparisons in the worst case. Since $Q$ is expected to be much smaller than $L$, 
the tree-search is more efficient than a simple search in the whole array which obviously requires $L-1$
comparisons.
The worst-case bottom-up update of a 3--level tree is depicted in Fig.~\ref{fig:ttree}.
Moreover, in the case of the present method, the kernel and therefore the 
residual coefficient update is localized in frequency as well. Therefore,
tournament-style trees are used for keeping track of maxima for individual
window time positions (across frequency bins) in a similar manner.

The trees provide the maximal coefficient $\cres_k(p_\text{max})$ of the 
coefficient-domain residual, where $p_\text{max}$ is given in terms of the dictionary 
$w_\text{max}$ and the position $m_\text{max}, n_\text{max}$.

\subsection{Fast Update Step With Real Atoms From Complex Multi-Gabor Dictionary}
When dealing with real atoms, pairs of conjugated atoms are involved.
The inner product of the currently selected atom with the conjugated partner
$\inprod{\atom_\text{max}}{\overline{\atom}_\text{max}}$
can actually be extracted from the kernel $\vect{h}_{w_\text{max},w_\text{max}}$.
Note that since the kernel has been truncated, we consider the inner product
to be zero if the conjugated atom is not in the range of the kernel update.
We consider the inner product to be zero also if the conjugated atom is missing.

Before the update step, $c_{\text{max}}$ is adjusted using \eqref{eq:dualatprod}.
The solution update step is performed as in Alg.~\ref{alg:CDmp} step \ref{item:solupdate}
while the error update step \ref{item:errupdate} uses \eqref{eq:projenfaster} in place of
the squared magnitude of the coefficient if the conjugated partner of $\atom_\text{max}$
is present.
The steps of the coefficient residual update in the coefficient domain are 
summarized in Alg.~\ref{alg:resupdate}.
Note the substraction of the kernel is performed for the conjugated atom as well, if necessary.

\begin{algorithm*}[th]
\SetKwFunction{resupdate}{singleDictionaryResidualUpdate}
\DontPrintSemicolon
\caption{ Approximate coefficient-domain MP residual update }
\label{alg:resupdate}
    \KwIn{ $c_\text{max}$, 
    $w_\text{max}$, $m_\text{max},
    n_\text{max}$, $\inprod{\atom_\text{max}}{\overline{\atom}_\text{max}}$,
    $a_1,\dots,a_W$, $M_1,\dots,M_W$,
    $\vect{h}_{w,v}$ $(w,v=1,\dots,W)$ }
    \KwOut{ Inner products to be updated $\coef^{\res_{k}}_1,\dots,\coef^{\res_{k}}_W$
    }

    \For{$w=1,\dots,W$}{
        $\coef^{\res_{k}}_{w}\leftarrow$\resupdate{$m_\mathrm{max}$,
        $w$, $\coef^{\res_{k}}_{w}$ }\;
        \If{$|\inprod{\atom_\mathrm{max}}{\overline{\atom}_\mathrm{max}}|>0$}{
            $m_\text{conj} = M_{w_\text{max}} - m_\text{max}$\;
            $\coef^{\res_{k}}_{w}\leftarrow$\resupdate{$m_\mathrm{conj}$, $w$, $\coef^{\res_{k}}_{w}$}\; }
    }
    \bigskip
  \SetKwProg{Pn}{Function}{:}{\KwRet}
  \Pn{$\coef_{w}\leftarrow$\resupdate{$m_\mathrm{max}$, $w$, $\coef_{w}$ }}{
        $a_\text{rat} \leftarrow   a_w/a_{w_\text{max}}$,
        $M_\text{rat} \leftarrow M_{w_\text{max}} / M_w$\;
        $a_\text{step} = a_\text{rat}$ or 1 if $a_\text{rat} < 1$,
        $M_\text{step} = M_\text{rat}$ or 1 if $M_\text{rat} < 1$\;
        \tcc{Determine index sets}
        Define $\mathcal{M}$, horizontal index set with stride 
        $M_\text{step}$ in cross kernel $\vect{h}_{w_\mathrm{max},w}$
        taking into the account the misalignment of the grids. \;
        Define $\mathcal{N}$, vertical index set in a similar way using stride
        $a_\text{step}$.\;
        Define $\mathcal{I}$, residual coefficient vector index set covered
        by the kernel.\;
    \tcc{Update the residual (as in \eqref{eq:fastkernupdate}) using truncated, subsampled and
    modulated cross-kernel
    \eqref{eq:hkj}:}
    $\coef_{w}(\mathcal{I}) = \coef_{w}(\mathcal{I}) -
    c_\mathrm{max}\vect{h}^{(m_\text{max})}_{w_\mathrm{max},w}
    (\mathcal{M},\mathcal{N})$\;
  }
\end{algorithm*}

\section{Comparison with MPTK}
\label{sec:mptk}
In order to showcase the efficiency of the proposed algorithm and its implementation,
in this section, we present a comparison with MPTK (version 0.7.0),
which is considered to be the fastest implementation available. To our knowledge there is no previous implementation of coefficient-domain MP that is able to decompose signals of the size considered here.
We measure the duration of the matching pursuit decomposition only.
From MPTK, we used the modified \texttt{mpd} utility tool.
The testing signal was the first channel from the file no. 39 from the SQAM database \cite{sqam}, which 
is a 137 seconds long piano recording sampled at 44.1~kHz totaling $6\cdot 10^6$ samples.
%
Both implementations link the same FFTW library \cite{fftw05} version 3.3 and were compiled
using the GCC (g++) compiler (version 7.2.0) with the \texttt{-Ofast} optimization flag
enabled. The creation of the FFTW plans and the computation of the kernels was
excluded from the measurements.
The specification of the PC the timing was performed on was
Intel\textregistered{} Core\texttrademark{} i5-4570 3.20~GHz, 16~GB RAM running 
Ubuntu 16.04.
The timing was done on an idle machine using the high-precision timer from the C++11 standard library
\texttt{chrono}. The data type was double precision floating point.
We used single and multi Gabor dictionaries with the Blackman window and various redundancies
$M_w/a_w$. The length of the window was always equal to $M$ (required by MPTK).
In the decomposition we performed $1.8\cdot10^5$ iterations. By fixing the number of iterations instead of a desired approximation estimate, we ensure that execution time is roughly independent of the considered signal. 
Table \ref{tab:timing} shows a comparison of execution times, in seconds, for a single Gabor dictionary,
numbers of bins  $M_w=512,\dots,8192$ (additionally also $16384$ for the proposed implementation)
and various hop sizes $a$ (and therefore redundancies). Additionally, a comparison of execution times using two multi-dictionaries is shown, 
each of which consists of five Gabor dictionaries (at redundancies $4$ and $8$ per dictionary).

\def\arraystretch{1.1}
\begin{table}[th]
\small
\centering
\begin{tabular}{|c|cccccc|}
\hline
 Bins ($M$) & $512$ & $1024$ & $2048$ & $4096$ & $8192$ & $16384$ \\
 \hline
\hline
  & \multicolumn{6}{|c|}{$a = M/4$} \\
 \hline
 MPTK & $3.96$ & $8.40$ & $17.0$ & $36.4$ & $75.7$ & -- \\
 Proposed & $0.92$ & $1.00$ & $1.02$ & $1.03$ & $1.03$ & $1.08$ \\
 \hline
  & \multicolumn{6}{|c|}{$a = M/8$} \\
 \hline
 MPTK & $6.73$ & $15.1$ & $30.2$ & $61.3$ & $147$ & -- \\
 Proposed & $1.70$ & $1.90$ & $1.95$ & $2.10$ & $2.20$ & $2.21$ \\
 \hline
  & \multicolumn{6}{|c|}{$a = M/16$} \\
 \hline
 MPTK & $13.2$ & $28.3$ & $56.8$ & $119$ & $274$ & -- \\
 Proposed & $3.20$ & $3.50$ & $4.00$ & $4.50$ & $4.60$ & $5.08$ \\
 \hline
  & \multicolumn{6}{|c|}{$a = M/32$} \\
 \hline
 MPTK & $23.2$ & $52.6$ & $110$ & $233$ & $530$ & -- \\
 Proposed & $6.35$ & $7.40$ & $7.90$ & $8.20$ & $9.60$ & $10.7$ \\
 \hline
 \hline
 Multi-Gabor & \multicolumn{3}{|c|}{$a = M/4$} & \multicolumn{3}{|c|}{$a = M/8$} \\
 \hline
 MPTK & \multicolumn{3}{|c|}{$142$} & \multicolumn{3}{|c|}{$285$} \\
 Proposed & \multicolumn{3}{|c|}{$10.6$} & \multicolumn{3}{|c|}{$20.9$} \\
 \hline
\end{tabular}
\caption{Execution time in seconds for MPTK and the proposed method on Gabor and Multi-Gabor dictionaries ($180$k selection steps). The Multi-Gabor dictionaries are a concatenation of the dictionaries with $M = 512,\ldots,8192$ at redundancies $a = M/4$ and $a = M/8$, respectively.}\label{tab:timing}
\vspace{-0.3cm}
\end{table}
\def\arraystretch{1.0}

In the tested setting, the proposed implementation clearly outperforms MPTK in 
terms of computational speed. The memory requirements are however notably higher since
the residual is stored in the coefficient domain and, additionally, the pre-computed 
kernels and the modulation factors must be stored as well. Hence, the proposed method requires
additional memory in an amount roughly proportional to the redundancy of the dictionary $\Dict$.
Note, however, that storage of the kernels and the modulation factors is independent of the signal 
length $L$ and thus increasingly insignificant.

\subsection{Convergence in practice}\label{sec:practconv}
\begin{figure}[th]
    \centering
    \includegraphics[width=6.8cm,height=4.5cm,trim={1.6cm 0 1.5cm 0},clip]{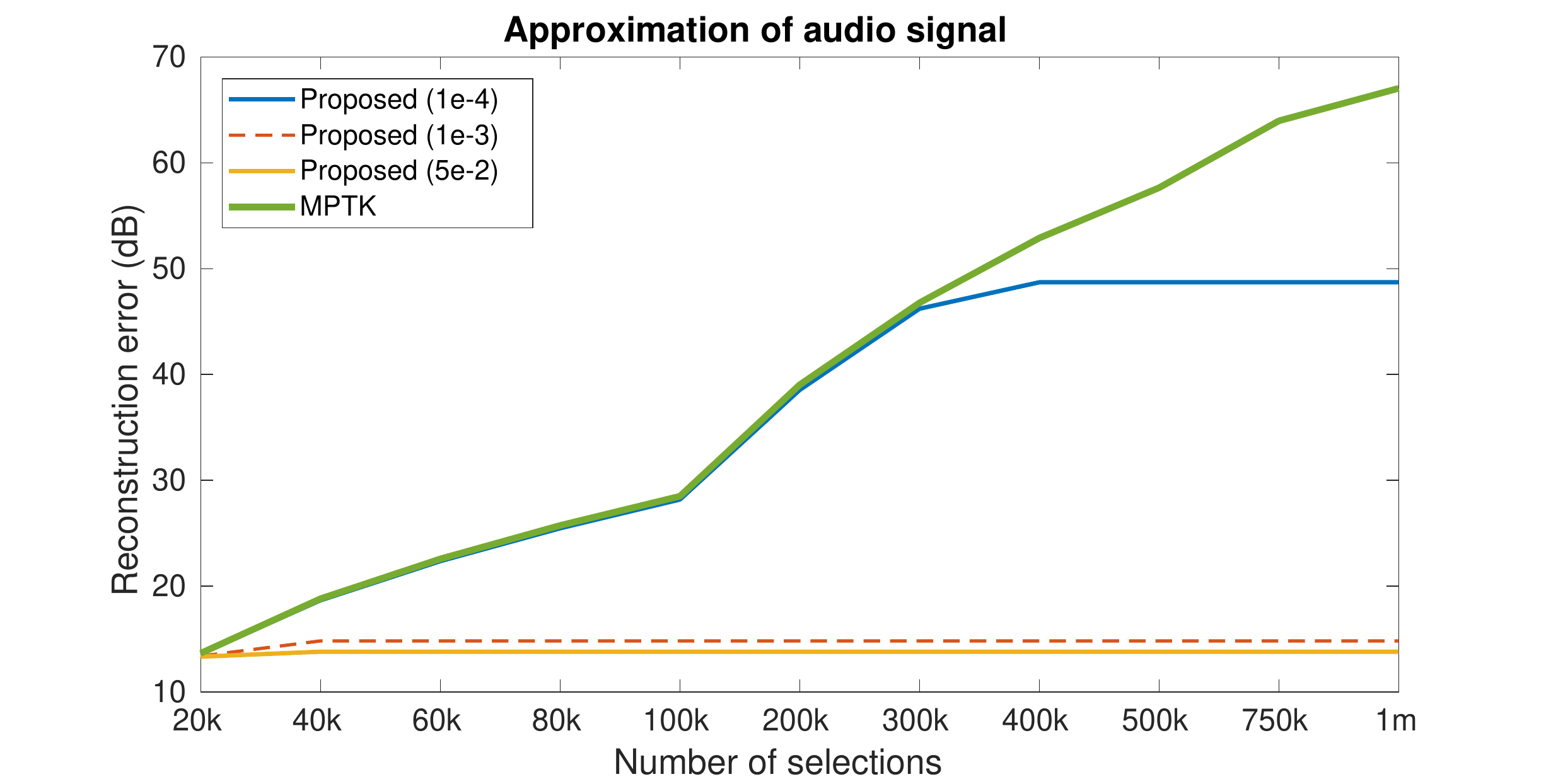}
    \includegraphics[width=6.8cm,height=4.5cm,trim={1cm 0 1.9cm 0},clip]{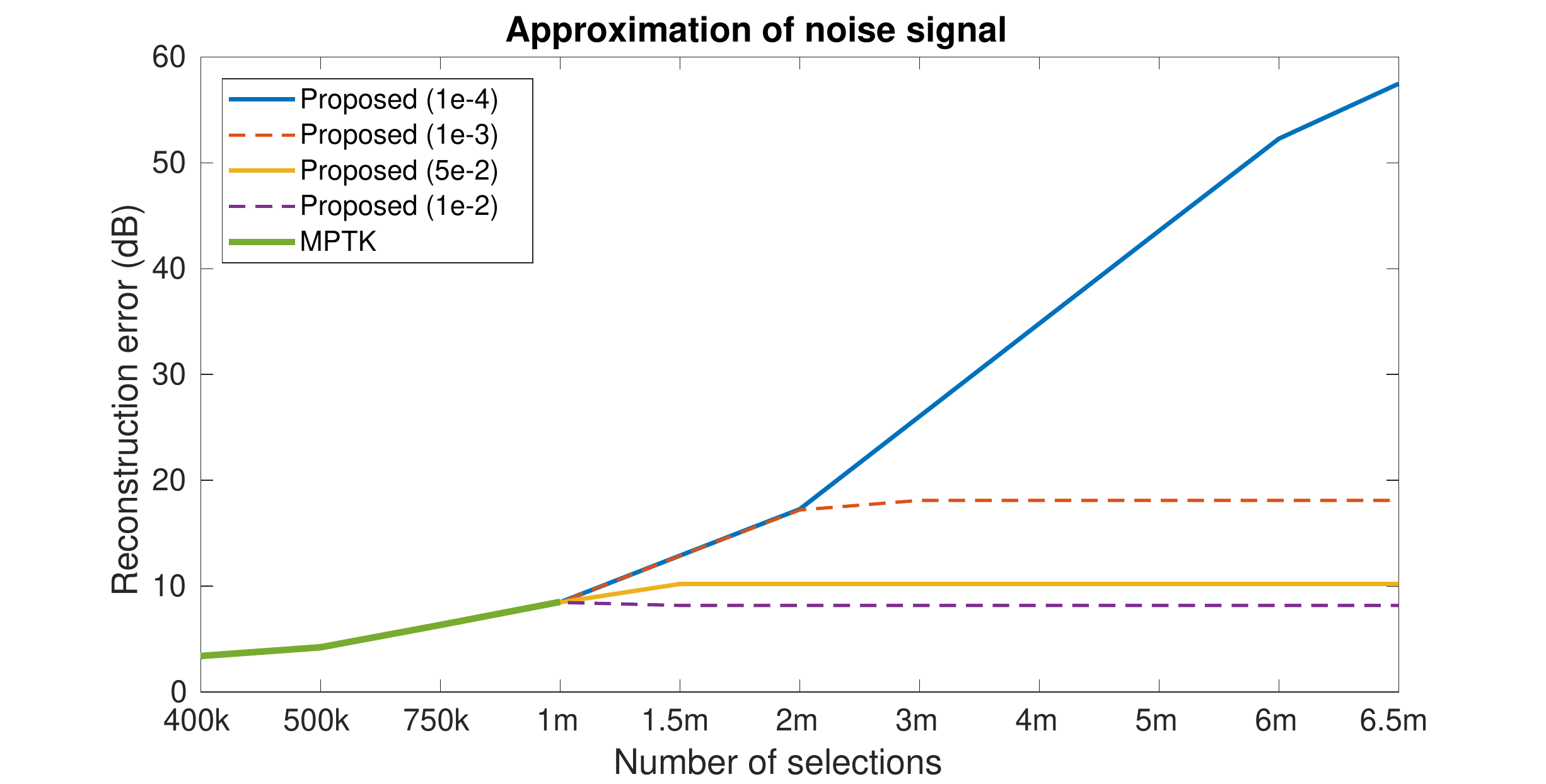}
   \caption{Approximation errors achieved by MPTK and the proposed method. The proposed method was evaluated for various threshold values $\epsilon$. For the audio signal (top) the proposed method terminated after $368235$ ($\epsilon = 10^{-4}$), $35115$ ($\epsilon = 10^{-3}$) and $29138$ ($\epsilon = 5\cdot10^{-3}$) steps, achieving an approximation quality of $48.72$, $14.84$ and $13.83$~dB, respectively. For the noise signal (bottom) the proposed method terminated after $6623689$ ($\epsilon = 10^{-4}$), $2160138$ ($\epsilon = 10^{-3}$), $1320287$ ($\epsilon = 5\cdot10^{-3}$) and $1240012$ ($\epsilon = 10^{-3}$) steps, achieving an approximation quality of $57.48$, $18.12$, $10.22$ and $8.19$~dB, respectively.}
    \label{fig:approxresults1}
    \vspace{-0.3cm}
  \end{figure}

  The results in Section \ref{ssec:approxupdate} suggest that the proposed approximate coefficient-domain update will not achieve arbitrarily small approximation error in general. Arbitrary approximation quality can only be achieved if resets are performed. To test the necessity of resets in practical applications, we compare the residual norm achieved after $k$ selection steps by MPTK and the proposed method. We do so for various values of $k$ and the truncation threshold $\epsilon$. The experiment was performed using a concatenation of $3$ Gabor dictionaries with Blackman window, $M_w = 512,1024,2048$ and $a_w = M_w/4$, as used in Section \ref{sec:mptk}. In addition to the audio test signal used previously, we also consider a pseudo-random Gaussian noise of equal length, generated in Matlab using the function \texttt{randn}. While the considered audio signal is approximately sparse and can be approximated with a moderate number of atoms, the noise signal is not sparse at all. Generally, our implementation of the proposed method is set to terminate when $\sum_{l=1}^k \cres_{l-1}(p_{l})^2 > \|\sig\|^2_2$. This condition simply means that the error estimate $E_k$ is negative and serves as a cheap indicator that further selection steps are expected to harm the approximation quality. Although this condition is not sufficient to prevent unproductive selection steps altogether, it was usually sufficient to prevent notable divergence of the matching pursuit estimate for various signals and threshold values, in particular for all presented experiments. The number of selection steps performed until this stopping condition is met depends, however, heavily on the considered signal, see Figure \ref{fig:approxresults1}. Instead, we observed similar final  approximation quality across different signals when this is the only stopping condition. We also state the final approximation error and number $k$ of selection steps performed before $\sum_{l=1}^k \cres_{l-1}(p_{l})^2 > \|\sig\|^2_2$. In all experiments, approximation quality of the proposed method follows MPTK closely, as long as $\sum_{l=1}^k \cres_{l-1}(p_{l})^2$ does not approach $\|\sig\|^2_2$. In practice, memory usage of MPTK grows linearly in the number of selection steps, such that we could only test MPTK approximation quality up to $10^6$ selection steps.
  
  When resets are activated, our Matlab interface relies on a reset criterion combining a maximum number of selections per reset and some relative error tolerance by default, see Section \ref{sec:module}. The default reset criterion depends on the kernel threshold and was selected to achieve desired approximation quality in reasonable time for audio signals and kernel thresholds in the range $\epsilon \in [10^{-2},10^{-4}]$. Therefore both conditions are set to conservative values and we do not claim optimality in terms of runtime or the number of selected atoms. 
  
  On the other hand, the number of steps before termination obtained in the test above, see Figure \ref{fig:approxresults1}, provides us with a better idea of how many \emph{good} selections we can perform before a reset is required. For the two tested signals, we can choose a less conservative reset condition. For a second experiment, we repeated the previous test, setting the reset criterion to a fixed number of iterations somewhat below the numbers obtained in the first test. With these values, we attempted once more to approximate both test signals. Specifically, resets were performed after $32\cdot 10^4,\ 3\cdot 10^4$ and $2\cdot 10^4$ selection steps for the kernel thresholds $10^{-4},\ 10^{-3}$ and $5\cdot10^{-3}$, respectively, on the audio test signal. For the noise signal, resets were performed after $7 \cdot 10^6,\ 2\cdot 10^6,\ 10^6$ and $8\cdot 10^5$ selection steps for thresholds $10^{-4},\ 10^{-3},\ 5\cdot10^{-3}$ and $10^{-3}$, respectively. In all tested cases, the error of the MP approximation follows the same curve, up to deviations of $\pm 1$~dB for the audio signal and $\pm 0.1$~dB for the noise signal (below $6\cdot 10^6$ selections). Up to $10^6$ selection steps (the tested range), approximation quality achieved by MPTK follows the same curve. Note that both signals are $6.04\cdot 10^6$ samples long, such that performing more than $6\cdot 10^6$ selection steps is unlikely to be practically relevant.
  
\section{Using the C module}
\label{sec:module}

Using the Large Time-Frequency Analysis Toolbox (LTFAT, \url{http://ltfat.github.io}), we provide the interface \texttt{multidgtrealmp.m} for using the proposed matching pursuit implementation from Matlab/GNU Octave. 

\begin{fundef}
  {\texttt{c = multidgtrealmp(f,dicts,errdb,maxit)}}\\
  
  Computes the MP decomposition \texttt{c} of input vector \texttt{f} with respect to multi-Gabor dictionary \texttt{dicts} given as cell array of triplets \texttt{\{g,a,M,$\ldots$\}} specifying the Gabor dictionary using the supported window type \texttt{g}, hop size \texttt{a} and \texttt{M} frequency channels. The window length will be set to \texttt{M}. The optional parameters \texttt{errdb} and \texttt{maxit} are used to specify the stopping criterion in terms of the targeted residual energy in dB and the maximum number of selection steps.
  
  The key-value pair \texttt{'kernthr',$\epsilon$} can be used to set the kernel threshold $\epsilon$. The default value is $10^{-4}$. 
  
  The flag \texttt{'reset'} activates resets. By default, the inner stopping criterion is derived from the kernel threshold $\epsilon$ and given by a maximum number of selections per reset and a tolerance for the relative error estimate $E_k/E^{\text{out}}$. Reset conditions can be set manually using the key-value pairs \texttt{'resetit', it} and \texttt{'reseterrdb', err}, where \texttt{it} is the maximum number of iterations per reset and the latter triggers a reset when $\texttt{err} > 10\log_{10}(E_k/E^{\text{out}})$.
  
  Further options are detailed in the documentation of \texttt{multidgtrealmp.m}.
\end{fundef}

\vspace{.1cm}
To demonstrate the usage of the \texttt{dgtrealmp} module of LIBLTFAT in C or C++ directly, we provide the example implementations \texttt{multigabormpd.cpp} 
at \url{https://github.com/ltfat/libltfat/tree/master/examples/multigabormp}. An overview of the module can be found at \url{http://ltfat.github.io/libltfat/group__multidgtrealmp}. The general workflow is as follows.

\begin{fundef}
  1) Initialize the parameter setup structure:\\
  {\footnotesize\texttt{int dgtrealmp\_parbuf\_init(dgtrealmp\_parbuf $^{**}$p)}}\\  
  
  \noindent 2) Add dictionaries to the setup structure using\\
  {\footnotesize\texttt{int dgtrealmp\_parbuf\_add\_firwin(dgtrealmp\_parbuf $^*$parbuf, LTFAT\_FIRWIN win, ltfat\_int gl, ltfat\_int a, ltfat\_int M)}}\\
  Here, \texttt{win} is a supported window type, \texttt{gl} is the window length and \texttt{a, M} are the hop size and number of frequency channels of the Gabor dictionary, respectively.\\
  
  \noindent 3) Optionally change MP parameters using the functions\\
  {\footnotesize\texttt{int dgtrealmp\_setparbuf\_$\ldots$(dgtrealmp\_parbuf $^*$parbuf, $\ldots$)}}\\
  By default, the algorithm targets an error estimate $10\log_{10}(E_k/\|\sig\|^2) \leq -40$ with at most $L/5$ selections at relative kernel threshold $10^{-4}$. Examples:\\
  {\footnotesize\texttt{int dgtrealmp\_setparbuf\_maxit(dgtrealmp\_parbuf $^*$parbuf, size\_t maxit)}}\\
  {\footnotesize\texttt{int dgtrealmp\_setparbuf\_errtoldb(dgtrealmp\_parbuf $^*$parbuf, double errtoldb)}}\\
  {\footnotesize\texttt{int dgtrealmp\_setparbuf\_kernrelthr (dgtrealmp\_parbuf $^*$parbuf, double thr)}}\\
  See \url{http://ltfat.github.io/libltfat/group__multidgtrealmp} for more options.\\
  
  \noindent 4) Initialize MP state for fixed input length \texttt{L}:\\
  {\footnotesize\texttt{int dgtrealmp\_init(dgtrealmp\_parbuf $^*$pb, ltfat\_int L, dgtrealmp\_state $^{**}$p)}}\\
  Precomputes the truncated kernels and cross-kernels, as well as all other data necessary to perform the MP decomposition.\\
  
  \noindent 5) Compute MP approximation for input vector \texttt{f}:\\
  {\footnotesize\texttt{int dgtrealmp\_execute(dgtrealmp\_state $^*$p, const LTFAT\_REAL f[], LTFAT\_COMPLEX $^*$cout[], LTFAT\_REAL fout[])}}\\
  The outputs \texttt{cout} and \texttt{fout} represent the final MP approximation in the coefficient and signal domains, respectively. This step may be repeated at will for different input vectors of the same length \texttt{L}.\\
  
  \noindent 6) Clean up:\\
  {\footnotesize\texttt{int dgtrealmp\_done(dgtrealmp\_state $^{**}$p)}}
\end{fundef}

\vspace{.1cm}
The function \texttt{dgtrealmp\_execute} computes the MP decomposition \texttt{cout} and synthesizes the signal domain approximation \texttt{fout}. The current implementation only supports resets via the Matlab/GNU Octave interface. In a future update, the reset procedure will be implemented directly into the C module.

\section{Conclusion And Outlook}
We have presented an accelerated MP algorithm alongside a reference implementation  
suitable for multi-Gabor dictionaries. 
Due to the structure of the Gram matrix of the multi-Gabor dictionary, 
the coefficient domain residual update step becomes very fast while
the memory requirements for storing the inner products between the atoms
remain constant with increasing signal length.
Moreover, the time and frequency locality of the residual update in turn allows faster
search for the maximum in the next iteration.
We have shown that the proposed method converges to the true solution, if a simple reset procedure is occasionally performed. 
Benchmarks show that, depending on the dictionary, our implementation is 3.5--70 times faster than 
the standard MPTK implementation. In the single dictionary case, the most notable feature is that
the execution time is virtually independent of the number of bins $M$ when the redundancy $M/a$ is fixed.
Moreover, as it turned out, MPTK could not handle dictionaries with the number of bins
higher than $8192$.
In our code, no explicit optimization techniques like exploiting the SIMD operations or
parallelization of the code were used, therefore it is possible that there is
still room for improvement.

Since the presented acceleration technique applies only to
the update step of the algorithm, it is also applicable to
various extensions of MP 
like the molecular MP \cite{da06}, perceptual MP \cite{lana08,chaneba14},
and guided MP \cite{za16} etc.
Extensions to the Local OMP \cite{magrbiva09,magrvabi11}, cyclic MP
\cite{stchr10,stchrgr11}, self-projected MP \cite{rerose17} and to the complementary MP \cite{ragu08,ragu10}
seem to be possible as well.
Future work may investigate the suitability of the proposed
implementation to sliding local MP \cite{da10} and real-time
MP mentioned in \cite{za16}.
We will also investigate the structure of the Gram matrix of a dictionary
consisting of an ensemble of  
real-valued windowed discrete cosine bases used by Ravelli~et.~al.~\cite{rarida08} to determine whether
similar acceleration technique is feasible.
%

\section*{Acknowledgment}
This work was supported by the Austrian Science Fund (FWF): Y\,551--N13 and I\,3067--N30.
The authors thank Bob L. Sturm for sharing his thoughts on the subject in 
a form of a blog \emph{Pursuits in the Null Space} and to the anonymous reviewers for valuable comments.



\bibliographystyle{ACM-Reference-Format}
\bibliography{project}


\begin{thebibliography}{47}


\ifx \showCODEN    \undefined \def \showCODEN     #1{\unskip}     \fi
\ifx \showDOI      \undefined \def \showDOI       #1{#1}\fi
\ifx \showISBNx    \undefined \def \showISBNx     #1{\unskip}     \fi
\ifx \showISBNxiii \undefined \def \showISBNxiii  #1{\unskip}     \fi
\ifx \showISSN     \undefined \def \showISSN      #1{\unskip}     \fi
\ifx \showLCCN     \undefined \def \showLCCN      #1{\unskip}     \fi
\ifx \shownote     \undefined \def \shownote      #1{#1}          \fi
\ifx \showarticletitle \undefined \def \showarticletitle #1{#1}   \fi
\ifx \showURL      \undefined \def \showURL       {\relax}        \fi
\providecommand\bibfield[2]{#2}
\providecommand\bibinfo[2]{#2}
\providecommand\natexlab[1]{#1}
\providecommand\showeprint[2][]{arXiv:#2}

\bibitem[\protect\citeauthoryear{Balazs}{Balazs}{2008}]%
        {xxlfinfram1}
\bibfield{author}{\bibinfo{person}{P. Balazs}.}
  \bibinfo{year}{2008}\natexlab{}.
\newblock \showarticletitle{Frames and finite dimensionality: {F}rame
  transformation, classification and algorithms}.
\newblock \bibinfo{journal}{\emph{Applied Mathematical Sciences}}
  \bibinfo{volume}{2}, \bibinfo{number}{41--44} (\bibinfo{year}{2008}),
  \bibinfo{pages}{2131--2144}.
\newblock


\bibitem[\protect\citeauthoryear{Bhattacharya and Depalle}{Bhattacharya and
  Depalle}{2014}]%
        {bhde14}
\bibfield{author}{\bibinfo{person}{G. Bhattacharya} {and} \bibinfo{person}{P.
  Depalle}.} \bibinfo{year}{2014}\natexlab{}.
\newblock \showarticletitle{Sparse denoising of audio by greedy time-frequency
  shrinkage}. In \bibinfo{booktitle}{\emph{Proc. IEEE Int. Conf. on Acoustics,
  Speech and Signal Processing (ICASSP)}}. \bibinfo{pages}{2898--2902}.
\newblock
\showISSN{1520-6149}


\bibitem[\protect\citeauthoryear{Blumensath and Davies}{Blumensath and
  Davies}{2008a}]%
        {blda08a}
\bibfield{author}{\bibinfo{person}{T. Blumensath} {and} \bibinfo{person}{M.~E.
  Davies}.} \bibinfo{year}{2008}\natexlab{a}.
\newblock \showarticletitle{Gradient pursuit for non-linear sparse signal
  modelling}. In \bibinfo{booktitle}{\emph{Proc. European Signal Processing
  Conference (EUSIPCO)}}. IEEE, \bibinfo{pages}{1--5}.
\newblock


\bibitem[\protect\citeauthoryear{Blumensath and Davies}{Blumensath and
  Davies}{2008b}]%
        {blda08b}
\bibfield{author}{\bibinfo{person}{T. Blumensath} {and} \bibinfo{person}{M.~E.
  Davies}.} \bibinfo{year}{2008}\natexlab{b}.
\newblock \showarticletitle{Gradient pursuits}.
\newblock \bibinfo{journal}{\emph{IEEE Tran. Signal Processing}}
  \bibinfo{volume}{56}, \bibinfo{number}{6} (\bibinfo{year}{2008}),
  \bibinfo{pages}{2370--2382}.
\newblock


\bibitem[\protect\citeauthoryear{Chardon, Necciari, and Balazs}{Chardon
  et~al\mbox{.}}{2014}]%
        {chaneba14}
\bibfield{author}{\bibinfo{person}{G. Chardon}, \bibinfo{person}{T. Necciari},
  {and} \bibinfo{person}{P. Balazs}.} \bibinfo{year}{2014}\natexlab{}.
\newblock \showarticletitle{Perceptual matching pursuit with {G}abor
  dictionaries and time-frequency masking}. In \bibinfo{booktitle}{\emph{Proc.
  IEEE Int. Conf. on Acoustics, Speech and Signal Processing (ICASSP)}}.
  \bibinfo{pages}{3102--3106}.
\newblock
\showISSN{1520-6149}


\bibitem[\protect\citeauthoryear{Daudet}{Daudet}{2006}]%
        {da06}
\bibfield{author}{\bibinfo{person}{L. Daudet}.}
  \bibinfo{year}{2006}\natexlab{}.
\newblock \showarticletitle{Sparse and structured decompositions of signals
  with the molecular matching pursuit}.
\newblock \bibinfo{journal}{\emph{IEEE Tran. Audio, Speech, and Language
  Processing}} \bibinfo{volume}{14}, \bibinfo{number}{5} (\bibinfo{date}{Sept}
  \bibinfo{year}{2006}), \bibinfo{pages}{1808--1816}.
\newblock
\showISSN{1558-7916}


\bibitem[\protect\citeauthoryear{Daudet}{Daudet}{2010}]%
        {da10}
\bibfield{author}{\bibinfo{person}{L. Daudet}.}
  \bibinfo{year}{2010}\natexlab{}.
\newblock \showarticletitle{Audio sparse decompositions in parallel: {L}et the
  greed be shared!}
\newblock \bibinfo{journal}{\emph{IEEE Signal Processing Magazine}}
  \bibinfo{volume}{27}, \bibinfo{number}{2} (\bibinfo{date}{March}
  \bibinfo{year}{2010}), \bibinfo{pages}{90--96}.
\newblock
\showISSN{1053-5888}


\bibitem[\protect\citeauthoryear{Davis, Mallat, and Avellaneda}{Davis
  et~al\mbox{.}}{1997}]%
        {damaav97}
\bibfield{author}{\bibinfo{person}{G. Davis}, \bibinfo{person}{S. Mallat},
  {and} \bibinfo{person}{M. Avellaneda}.} \bibinfo{year}{1997}\natexlab{}.
\newblock \showarticletitle{Adaptive greedy approximations}.
\newblock \bibinfo{journal}{\emph{Constructive Approximation}}
  \bibinfo{volume}{13}, \bibinfo{number}{1} (\bibinfo{year}{1997}),
  \bibinfo{pages}{57--98}.
\newblock
\showISSN{1432-0940}


\bibitem[\protect\citeauthoryear{Davis, Mallat, and Zhang}{Davis
  et~al\mbox{.}}{1994}]%
        {damazh94}
\bibfield{author}{\bibinfo{person}{G.~M. Davis}, \bibinfo{person}{S.~G.
  Mallat}, {and} \bibinfo{person}{Z. Zhang}.} \bibinfo{year}{1994}\natexlab{}.
\newblock \showarticletitle{Adaptive time-frequency decompositions}.
\newblock \bibinfo{journal}{\emph{Optical Engineering}} \bibinfo{volume}{33},
  \bibinfo{number}{7} (\bibinfo{year}{1994}), \bibinfo{pages}{2183--2191}.
\newblock
\showISBNx{0091-3286}


\bibitem[\protect\citeauthoryear{Derrien}{Derrien}{2007}]%
        {de07}
\bibfield{author}{\bibinfo{person}{O. Derrien}.}
  \bibinfo{year}{2007}\natexlab{}.
\newblock \showarticletitle{{Time-scaling of audio signals with muti-scale
  Gabor analysis}}. In \bibinfo{booktitle}{\emph{Proc. Int. Conf. Digital Audio
  Effects (DAFx--07)}}. \bibinfo{address}{Bordeaux, France}.
\newblock


\bibitem[\protect\citeauthoryear{DeVore and Temlyakov}{DeVore and
  Temlyakov}{1996}]%
        {dete96}
\bibfield{author}{\bibinfo{person}{R.~A. DeVore} {and} \bibinfo{person}{V.~N.
  Temlyakov}.} \bibinfo{year}{1996}\natexlab{}.
\newblock \showarticletitle{Some remarks on greedy algorithms}.
\newblock \bibinfo{journal}{\emph{Advances in Computational Mathematics}}
  \bibinfo{volume}{5}, \bibinfo{number}{1} (\bibinfo{date}{01 Dec}
  \bibinfo{year}{1996}), \bibinfo{pages}{173--187}.
\newblock
\showISSN{1572-9044}


\bibitem[\protect\citeauthoryear{Durka}{Durka}{2007}]%
        {du07}
\bibfield{author}{\bibinfo{person}{P. Durka}.} \bibinfo{year}{2007}\natexlab{}.
\newblock \bibinfo{booktitle}{\emph{Matching pursuit and unification in {EEG}
  analysis}}.
\newblock \bibinfo{publisher}{Artech House, Inc.} 184 pages.
\newblock
\showISBNx{978-1-58053-304-1}


\bibitem[\protect\citeauthoryear{Ferrando, Kolasa, and
  Kova\v{c}evi\'{c}}{Ferrando et~al\mbox{.}}{2002}]%
        {fe02}
\bibfield{author}{\bibinfo{person}{S.~E. Ferrando}, \bibinfo{person}{L.~A.
  Kolasa}, {and} \bibinfo{person}{N. Kova\v{c}evi\'{c}}.}
  \bibinfo{year}{2002}\natexlab{}.
\newblock \showarticletitle{Algorithm 820: A flexible implementation of
  matching pursuit for {G}abor functions on the interval}.
\newblock \bibinfo{journal}{\emph{ACM Trans. Math. Softw.}}
  \bibinfo{volume}{28}, \bibinfo{number}{3} (\bibinfo{date}{Sept.}
  \bibinfo{year}{2002}), \bibinfo{pages}{337--353}.
\newblock
\showISSN{0098-3500}


\bibitem[\protect\citeauthoryear{Frigo and Johnson}{Frigo and Johnson}{2005}]%
        {fftw05}
\bibfield{author}{\bibinfo{person}{M. Frigo} {and} \bibinfo{person}{S.~G.
  Johnson}.} \bibinfo{year}{2005}\natexlab{}.
\newblock \showarticletitle{The design and implementation of {FFTW3}}.
\newblock \bibinfo{journal}{\emph{Proc. of the IEEE}} \bibinfo{volume}{93},
  \bibinfo{number}{2} (\bibinfo{year}{2005}), \bibinfo{pages}{216--231}.
\newblock
\newblock
\shownote{Special issue on ``Program Generation, Optimization, and Platform
  Adaptation''.}


\bibitem[\protect\citeauthoryear{Gribonval}{Gribonval}{2001}]%
        {gr01}
\bibfield{author}{\bibinfo{person}{R. Gribonval}.}
  \bibinfo{year}{2001}\natexlab{}.
\newblock \showarticletitle{Fast matching pursuit with a multiscale dictionary
  of {G}aussian chirps}.
\newblock \bibinfo{journal}{\emph{IEEE Tran. Signal Processing}}
  \bibinfo{volume}{49}, \bibinfo{number}{5} (\bibinfo{date}{May}
  \bibinfo{year}{2001}), \bibinfo{pages}{994--1001}.
\newblock
\showISSN{1053-587X}


\bibitem[\protect\citeauthoryear{Gribonval}{Gribonval}{2002}]%
        {gr02}
\bibfield{author}{\bibinfo{person}{R. Gribonval}.}
  \bibinfo{year}{2002}\natexlab{}.
\newblock \showarticletitle{Sparse decomposition of stereo signals with
  matching pursuit and application to blind separation of more than two sources
  from a stereo mixture}. In \bibinfo{booktitle}{\emph{IEEE Int. Conf.
  Acoustics, Speech, and Signal Processing}}, Vol.~\bibinfo{volume}{3}.
  \bibinfo{pages}{III--3057--III--3060}.
\newblock
\showISSN{1520-6149}


\bibitem[\protect\citeauthoryear{Gribonval, Depalle, Rodet, Bacry, and
  Mallat}{Gribonval et~al\mbox{.}}{1996}]%
        {grderobama96}
\bibfield{author}{\bibinfo{person}{R. Gribonval}, \bibinfo{person}{P. Depalle},
  \bibinfo{person}{X. Rodet}, \bibinfo{person}{E. Bacry}, {and}
  \bibinfo{person}{S. Mallat}.} \bibinfo{year}{1996}\natexlab{}.
\newblock \showarticletitle{{Sound signals decomposition using a high
  resolution matching pursuit}}. In \bibinfo{booktitle}{\emph{{Proc. Int.
  Computer Music Conf. (ICMC'96)}}}. \bibinfo{pages}{293--296}.
\newblock


\bibitem[\protect\citeauthoryear{Gribonval, {Figueras i Ventura}, and
  Vandergheynst}{Gribonval et~al\mbox{.}}{2006}]%
        {grfiva06}
\bibfield{author}{\bibinfo{person}{R. Gribonval}, \bibinfo{person}{R.~M.
  {Figueras i Ventura}}, {and} \bibinfo{person}{P. Vandergheynst}.}
  \bibinfo{year}{2006}\natexlab{}.
\newblock \showarticletitle{A simple test to check the optimality of a sparse
  signal approximation}.
\newblock \bibinfo{journal}{\emph{Signal Processing}} \bibinfo{volume}{86},
  \bibinfo{number}{3} (\bibinfo{year}{2006}), \bibinfo{pages}{496 -- 510}.
\newblock
\showISSN{0165-1684}


\bibitem[\protect\citeauthoryear{Gribonval and Vandergheynst}{Gribonval and
  Vandergheynst}{2006}]%
        {grva06}
\bibfield{author}{\bibinfo{person}{R. Gribonval} {and} \bibinfo{person}{P.
  Vandergheynst}.} \bibinfo{year}{2006}\natexlab{}.
\newblock \showarticletitle{On the exponential convergence of matching pursuits
  in quasi-incoherent dictionaries}.
\newblock \bibinfo{journal}{\emph{IEEE Tran. Information Theory}}
  \bibinfo{volume}{52}, \bibinfo{number}{1} (\bibinfo{date}{Jan}
  \bibinfo{year}{2006}), \bibinfo{pages}{255--261}.
\newblock
\showISSN{0018-9448}


\bibitem[\protect\citeauthoryear{Gr\"{o}chenig}{Gr\"{o}chenig}{2001}]%
        {chabible}
\bibfield{author}{\bibinfo{person}{K. Gr\"{o}chenig}.}
  \bibinfo{year}{2001}\natexlab{}.
\newblock \bibinfo{booktitle}{\emph{Foundations of time-frequency analysis}}.
\newblock \bibinfo{publisher}{Birkh\"{a}user}, \bibinfo{address}{Boston, Basel,
  Berlin}.
\newblock
\showISBNx{0-8176-4022-3}


\bibitem[\protect\citeauthoryear{Krstulovi\'c and Gribonval}{Krstulovi\'c and
  Gribonval}{2006}]%
        {krgr06}
\bibfield{author}{\bibinfo{person}{S. Krstulovi\'c} {and} \bibinfo{person}{R.
  Gribonval}.} \bibinfo{year}{2006}\natexlab{}.
\newblock \showarticletitle{{MPTK}: Matching pursuit made tractable}. In
  \bibinfo{booktitle}{\emph{Proc. Int. Conf. on Acoustics Speech and Signal
  Processing ICASSP 2006}}, Vol.~\bibinfo{volume}{3}.
  \bibinfo{pages}{III--496--III--499}.
\newblock
\showISSN{1520-6149}


\bibitem[\protect\citeauthoryear{Lahdili, Najaf-Zadeh, Pichevar, and
  Thibault}{Lahdili et~al\mbox{.}}{2008}]%
        {lana08}
\bibfield{author}{\bibinfo{person}{H. Lahdili}, \bibinfo{person}{H.
  Najaf-Zadeh}, \bibinfo{person}{R. Pichevar}, {and} \bibinfo{person}{L.
  Thibault}.} \bibinfo{year}{2008}\natexlab{}.
\newblock \showarticletitle{Perceptual matching pursuit for audio coding}. In
  \bibinfo{booktitle}{\emph{Audio Engineering Society Convention 124}}.
\newblock


\bibitem[\protect\citeauthoryear{{Le Roux}, Kameoka, Ono, and Sagayama}{{Le
  Roux} et~al\mbox{.}}{2010}]%
        {leroux10}
\bibfield{author}{\bibinfo{person}{J. {Le Roux}}, \bibinfo{person}{H. Kameoka},
  \bibinfo{person}{N. Ono}, {and} \bibinfo{person}{S. Sagayama}.}
  \bibinfo{year}{2010}\natexlab{}.
\newblock \showarticletitle{Fast signal reconstruction from magnitude {STFT}
  spectrogram based on spectrogram consistency}. In
  \bibinfo{booktitle}{\emph{Proc. 13th Int. Conf. on Digital Audio Effects
  (DAFx-10)}}. \bibinfo{pages}{397--403}.
\newblock


\bibitem[\protect\citeauthoryear{Leveau and Daudet}{Leveau and Daudet}{2006}]%
        {leda06}
\bibfield{author}{\bibinfo{person}{P. Leveau} {and} \bibinfo{person}{L.
  Daudet}.} \bibinfo{year}{2006}\natexlab{}.
\newblock \showarticletitle{Multi-resolution partial tracking with modified
  matching pursuit}. In \bibinfo{booktitle}{\emph{Proc. 14th European Signal
  Processing Conference}}. \bibinfo{pages}{1--4}.
\newblock
\showISSN{2219-5491}


\bibitem[\protect\citeauthoryear{Mailh\'{e}, Gribonval, Bimbot, and
  Vandergheynst}{Mailh\'{e} et~al\mbox{.}}{2009}]%
        {magrbiva09}
\bibfield{author}{\bibinfo{person}{B. Mailh\'{e}}, \bibinfo{person}{R.
  Gribonval}, \bibinfo{person}{F. Bimbot}, {and} \bibinfo{person}{P.
  Vandergheynst}.} \bibinfo{year}{2009}\natexlab{}.
\newblock \showarticletitle{A low complexity orthogonal matching pursuit for
  sparse signal approximation with shift-invariant dictionaries}. In
  \bibinfo{booktitle}{\emph{Proc. IEEE Int. Conf. on Acoustics, Speech and
  Signal Processing (ICAASP)}}. \bibinfo{pages}{3445--3448}.
\newblock
\showISSN{1520-6149}


\bibitem[\protect\citeauthoryear{Mailh\'{e}, Gribonval, Vandergheynst, and
  Bimbot}{Mailh\'{e} et~al\mbox{.}}{2011}]%
        {magrvabi11}
\bibfield{author}{\bibinfo{person}{B. Mailh\'{e}}, \bibinfo{person}{R.
  Gribonval}, \bibinfo{person}{P. Vandergheynst}, {and} \bibinfo{person}{F.
  Bimbot}.} \bibinfo{year}{2011}\natexlab{}.
\newblock \showarticletitle{Fast orthogonal sparse approximation algorithms
  over local dictionaries}.
\newblock \bibinfo{journal}{\emph{Signal Processing}} \bibinfo{volume}{91},
  \bibinfo{number}{12} (\bibinfo{year}{2011}), \bibinfo{pages}{2822 -- 2835}.
\newblock
\showISSN{0165-1684}


\bibitem[\protect\citeauthoryear{Mallat}{Mallat}{2008}]%
        {wtour}
\bibfield{author}{\bibinfo{person}{Stphane Mallat}.}
  \bibinfo{year}{2008}\natexlab{}.
\newblock \bibinfo{booktitle}{\emph{A Wavelet Tour of Signal Processing, Third
  Edition: {T}he Sparse Way} (\bibinfo{edition}{3rd} ed.)}.
\newblock \bibinfo{publisher}{Academic Press}.
\newblock
\showISBNx{0123743702, 9780123743701}


\bibitem[\protect\citeauthoryear{Mallat and Zhang}{Mallat and Zhang}{1993}]%
        {mazh93}
\bibfield{author}{\bibinfo{person}{S.~G. Mallat} {and} \bibinfo{person}{Z.
  Zhang}.} \bibinfo{year}{1993}\natexlab{}.
\newblock \showarticletitle{Matching pursuits with time-frequency
  dictionaries}.
\newblock \bibinfo{journal}{\emph{IEEE Tran. Signal Processing}}
  \bibinfo{volume}{41}, \bibinfo{number}{12} (\bibinfo{date}{Dec}
  \bibinfo{year}{1993}), \bibinfo{pages}{3397--3415}.
\newblock
\showISSN{1053-587X}


\bibitem[\protect\citeauthoryear{Pati, Rezaiifar, and Krishnaprasad}{Pati
  et~al\mbox{.}}{1993}]%
        {parekr93}
\bibfield{author}{\bibinfo{person}{Y.~C. Pati}, \bibinfo{person}{R. Rezaiifar},
  {and} \bibinfo{person}{P.~S. Krishnaprasad}.}
  \bibinfo{year}{1993}\natexlab{}.
\newblock \showarticletitle{Orthogonal matching pursuit: {R}ecursive function
  approximation with applications to wavelet decomposition}. In
  \bibinfo{booktitle}{\emph{Proc. 27th Asilomar Conference on Signals, Systems
  and Computers}}. \bibinfo{pages}{40--44 vol.1}.
\newblock
\showISSN{1058-6393}


\bibitem[\protect\citeauthoryear{Plumbley, Blumensath, Daudet, Gribonval, and
  Davies}{Plumbley et~al\mbox{.}}{2010}]%
        {plbldagrda10}
\bibfield{author}{\bibinfo{person}{M.~D. Plumbley}, \bibinfo{person}{T.
  Blumensath}, \bibinfo{person}{L. Daudet}, \bibinfo{person}{R. Gribonval},
  {and} \bibinfo{person}{M.~E. Davies}.} \bibinfo{year}{2010}\natexlab{}.
\newblock \showarticletitle{Sparse representations in audio and music: From
  coding to source separation}.
\newblock \bibinfo{journal}{\emph{Proc. IEEE}} \bibinfo{volume}{98},
  \bibinfo{number}{6} (\bibinfo{date}{June} \bibinfo{year}{2010}),
  \bibinfo{pages}{995--1005}.
\newblock
\showISSN{0018-9219}


\bibitem[\protect\citeauthoryear{Portnoff}{Portnoff}{1976}]%
        {po76}
\bibfield{author}{\bibinfo{person}{M.~R. Portnoff}.}
  \bibinfo{year}{1976}\natexlab{}.
\newblock \showarticletitle{Implementation of the digital phase vocoder using
  the fast Fourier transform}.
\newblock \bibinfo{journal}{\emph{IEEE Tran. Acoustics, Speech, and Signal
  Processing}} \bibinfo{volume}{24}, \bibinfo{number}{3} (\bibinfo{date}{Jun}
  \bibinfo{year}{1976}), \bibinfo{pages}{243--248}.
\newblock
\showISSN{0096-3518}


\bibitem[\protect\citeauthoryear{Pr\r{u}\v{s}a, S{\o}ndergaard, Holighaus,
  Wiesmeyr, and Balazs}{Pr\r{u}\v{s}a et~al\mbox{.}}{2014}]%
        {ltfatnote030}
\bibfield{author}{\bibinfo{person}{Z. Pr\r{u}\v{s}a}, \bibinfo{person}{P.~L.
  S{\o}ndergaard}, \bibinfo{person}{N. Holighaus}, \bibinfo{person}{Ch.
  Wiesmeyr}, {and} \bibinfo{person}{P. Balazs}.}
  \bibinfo{year}{2014}\natexlab{}.
\newblock \showarticletitle{{The large time-frequency analysis toolbox 2.0}}.
\newblock In \bibinfo{booktitle}{\emph{Sound, Music, and Motion}}.
  \bibinfo{publisher}{Springer International Publishing},
  \bibinfo{pages}{419--442}.
\newblock
\showISBNx{978-3-319-12975-4}


\bibitem[\protect\citeauthoryear{Rath and Guillemot}{Rath and
  Guillemot}{2008}]%
        {ragu08}
\bibfield{author}{\bibinfo{person}{G. Rath} {and} \bibinfo{person}{C.
  Guillemot}.} \bibinfo{year}{2008}\natexlab{}.
\newblock \showarticletitle{A complementary matching pursuit algorithm for
  sparse approximation}. In \bibinfo{booktitle}{\emph{Proc. 16th European
  Signal Processing Conference (EUSIPCO)}}. \bibinfo{pages}{1--5}.
\newblock
\showISSN{2219-5491}


\bibitem[\protect\citeauthoryear{Rath and Guillemot}{Rath and
  Guillemot}{2010}]%
        {ragu10}
\bibfield{author}{\bibinfo{person}{G. Rath} {and} \bibinfo{person}{Ch.
  Guillemot}.} \bibinfo{year}{2010}\natexlab{}.
\newblock \showarticletitle{On a simple derivation of the complementary
  matching pursuit}.
\newblock \bibinfo{journal}{\emph{Signal Processing}} \bibinfo{volume}{90},
  \bibinfo{number}{2} (\bibinfo{year}{2010}), \bibinfo{pages}{702 -- 706}.
\newblock
\showISSN{0165-1684}


\bibitem[\protect\citeauthoryear{Ravelli, Richard, and Daudet}{Ravelli
  et~al\mbox{.}}{2008}]%
        {rarida08}
\bibfield{author}{\bibinfo{person}{E. Ravelli}, \bibinfo{person}{G. Richard},
  {and} \bibinfo{person}{L. Daudet}.} \bibinfo{year}{2008}\natexlab{}.
\newblock \showarticletitle{Union of {MDCT} bases for audio coding}.
\newblock \bibinfo{journal}{\emph{IEEE Tran. Audio, Speech, and Language
  Processing}} \bibinfo{volume}{16}, \bibinfo{number}{8} (\bibinfo{date}{Nov}
  \bibinfo{year}{2008}), \bibinfo{pages}{1361--1372}.
\newblock
\showISSN{1558-7916}


\bibitem[\protect\citeauthoryear{Rebollo{-}Neira, Rozlo\v{z}n\'{i}k, and
  Sasmal}{Rebollo{-}Neira et~al\mbox{.}}{2017}]%
        {rerose17}
\bibfield{author}{\bibinfo{person}{L. Rebollo{-}Neira}, \bibinfo{person}{M.
  Rozlo\v{z}n\'{i}k}, {and} \bibinfo{person}{P. Sasmal}.}
  \bibinfo{year}{2017}\natexlab{}.
\newblock \showarticletitle{Analysis of a low memory implementation of the
  orthogonal matching pursuit greedy strategy}.
\newblock \bibinfo{journal}{\emph{CoRR}}  \bibinfo{volume}{abs/1609.00053v2}
  (\bibinfo{year}{2017}).
\newblock


\bibitem[\protect\citeauthoryear{Rish and Grabarnik}{Rish and
  Grabarnik}{2015}]%
        {rish15}
\bibfield{author}{\bibinfo{person}{I. Rish} {and} \bibinfo{person}{G.
  Grabarnik}.} \bibinfo{year}{2015}\natexlab{}.
\newblock \bibinfo{booktitle}{\emph{Sparse modeling: {T}heory, algorithms, and
  applications}}.
\newblock \bibinfo{publisher}{CRC Press}. 253 pages.
\newblock
\showISBNx{978-1-4398-2870-0}


\bibitem[\protect\citeauthoryear{S{\o}ndergaard}{S{\o}ndergaard}{2012}]%
        {ltfatnote011}
\bibfield{author}{\bibinfo{person}{P.~L. S{\o}ndergaard}.}
  \bibinfo{year}{2012}\natexlab{}.
\newblock \showarticletitle{Efficient algorithms for the discrete {G}abor
  transform with a long {FIR} window}.
\newblock \bibinfo{journal}{\emph{J.\ Fourier Anal.\ Appl.}}
  \bibinfo{volume}{18}, \bibinfo{number}{3} (\bibinfo{year}{2012}),
  \bibinfo{pages}{456--470}.
\newblock


\bibitem[\protect\citeauthoryear{S{\o}ndergaard, Torr\'esani, and
  Balazs}{S{\o}ndergaard et~al\mbox{.}}{2012}]%
        {ltfatnote015}
\bibfield{author}{\bibinfo{person}{P.~L. S{\o}ndergaard}, \bibinfo{person}{B.
  Torr\'esani}, {and} \bibinfo{person}{P. Balazs}.}
  \bibinfo{year}{2012}\natexlab{}.
\newblock \showarticletitle{{The linear time frequency analysis toolbox}}.
\newblock \bibinfo{journal}{\emph{International Journal of Wavelets,
  Multiresolution Analysis and Information Processing}} \bibinfo{volume}{10},
  \bibinfo{number}{4} (\bibinfo{year}{2012}), \bibinfo{pages}{1--27}.
\newblock


\bibitem[\protect\citeauthoryear{Sturm}{Sturm}{2009}]%
        {st09}
\bibfield{author}{\bibinfo{person}{B.~L. Sturm}.}
  \bibinfo{year}{2009}\natexlab{}.
\newblock \emph{\bibinfo{title}{Sparse approximation and atomic decomposition:
  {C}onsidering atom Interactions in evaluating and building signal
  representations}}.
\newblock \bibinfo{thesistype}{Ph.D. Dissertation}. \bibinfo{school}{University
  of California}.
\newblock


\bibitem[\protect\citeauthoryear{Sturm and Christensen}{Sturm and
  Christensen}{2010}]%
        {stchr10}
\bibfield{author}{\bibinfo{person}{B.~L. Sturm} {and} \bibinfo{person}{M.~G.
  Christensen}.} \bibinfo{year}{2010}\natexlab{}.
\newblock \showarticletitle{Cyclic matching pursuits with multiscale
  time-frequency dictionaries}. In \bibinfo{booktitle}{\emph{Conf. Record of
  the 44th Asilomar Conference on Signals, Systems and Computers}}.
  \bibinfo{pages}{581--585}.
\newblock
\showISSN{1058-6393}


\bibitem[\protect\citeauthoryear{Sturm, Christensen, and Gribonval}{Sturm
  et~al\mbox{.}}{2011}]%
        {stchrgr11}
\bibfield{author}{\bibinfo{person}{B.~L. Sturm}, \bibinfo{person}{M.~G.
  Christensen}, {and} \bibinfo{person}{R. Gribonval}.}
  \bibinfo{year}{2011}\natexlab{}.
\newblock \showarticletitle{Cyclic pure greedy algorithms for recovering
  compressively sampled sparse signals}. In
  \bibinfo{booktitle}{\emph{Conference Record of the 45th Asilomar Conference
  on Signals, Systems and Computers}}. IEEE, \bibinfo{pages}{1143--1147}.
\newblock


\bibitem[\protect\citeauthoryear{Sturm, Daudet, and Roads}{Sturm
  et~al\mbox{.}}{2006}]%
        {stdacu06}
\bibfield{author}{\bibinfo{person}{B.~L. Sturm}, \bibinfo{person}{L. Daudet},
  {and} \bibinfo{person}{C. Roads}.} \bibinfo{year}{2006}\natexlab{}.
\newblock \showarticletitle{Pitch-shifting audio signals using sparse atomic
  approximations}. In \bibinfo{booktitle}{\emph{Proc. 1st ACM Workshop on Audio
  and Music Computing Multimedia}} (Santa Barbara, California, USA)
  \emph{(\bibinfo{series}{AMCMM '06})}. \bibinfo{publisher}{ACM},
  \bibinfo{address}{New York, NY, USA}, \bibinfo{pages}{45--52}.
\newblock
\showISBNx{1-59593-501-0}


\bibitem[\protect\citeauthoryear{Sturm and Gibson}{Sturm and Gibson}{2006}]%
        {stgi06}
\bibfield{author}{\bibinfo{person}{B.~L. Sturm} {and} \bibinfo{person}{J.~D.
  Gibson}.} \bibinfo{year}{2006}\natexlab{}.
\newblock \showarticletitle{Matching pursuit decompositions of non-noisy speech
  signals using several dictionaries}. In \bibinfo{booktitle}{\emph{Proc. IEEE
  Int. Conf. on Acoustics Speech and Signal Processing}},
  Vol.~\bibinfo{volume}{3}. \bibinfo{pages}{III--III}.
\newblock
\showISSN{1520-6149}


\bibitem[\protect\citeauthoryear{{The European Broadcasting Union}}{{The
  European Broadcasting Union}}{2008}]%
        {sqam}
\bibfield{author}{\bibinfo{person}{{The European Broadcasting Union}}.}
  \bibinfo{year}{2008}\natexlab{}.
\newblock \bibinfo{booktitle}{\emph{Tech 3253: {Sound quality assessment
  material} recordings for subjective tests}}.
\newblock \bibinfo{type}{{T}echnical {R}eport}. \bibinfo{address}{Geneva}.
\newblock
\urldef\tempurl%
\url{https://tech.ebu.ch/docs/tech/tech3253.pdf}
\showURL{%
\tempurl}


\bibitem[\protect\citeauthoryear{Yaghoobi and Davies}{Yaghoobi and
  Davies}{2009}]%
        {yada09}
\bibfield{author}{\bibinfo{person}{M. Yaghoobi} {and} \bibinfo{person}{M.~E.
  Davies}.} \bibinfo{year}{2009}\natexlab{}.
\newblock \bibinfo{booktitle}{\emph{Fast and scalable: A survey on sparse
  approximation methods}}.
\newblock \bibinfo{type}{{T}echnical {R}eport}. \bibinfo{institution}{The
  University of Edinburgh}.
\newblock


\bibitem[\protect\citeauthoryear{Zantalis}{Zantalis}{2016}]%
        {za16}
\bibfield{author}{\bibinfo{person}{D. Zantalis}.}
  \bibinfo{year}{2016}\natexlab{}.
\newblock \emph{\bibinfo{title}{Guided matching pursuit and its application to
  sound source separation}}.
\newblock \bibinfo{thesistype}{Ph.D. Dissertation}. \bibinfo{school}{University
  of York}.
\newblock


\end{thebibliography}

\appendix

\section{Proofs}\label{app:proof1}

\subsection{Proof of Theorem \ref{thm:Resdecrease}}
We begin with some preparation. Recall that $\res_k = \sig - \Dict \coef_k$. By linearity, 
$\res_{k+1} = \res_k - \Dict (\coef_{k+1}-\coef_k) = \res_k - \cres_k(p_{k+1})\atom_{p_{k+1}}$.

Further, we have $\Dict^\ast \res_k = \cres_0 - \Gram \coef_k$ and 
$\cres_k = \cres_0 - \Gram_\epsilon \coef_k$. With $\sumcres_k := \sum_{l=1}^{k} |\cres_{l-1}(p_{l})|$, it is easily verified that 
$|(\Gram - \Gram_\epsilon) \coef_k|(p) \leq \epsilon \sumcres_k$, for all $p$, and consequently, 
$|\langle \res_k, \atom_p\rangle - \cres_k(p)| \leq \epsilon \sumcres_k$ as well.

By assumption, $\epsilon \sumcres_k \leq \delta |\cres_k(p_{k+1})|$, with $0<\delta< 1/2$. We have 
\begin{equation}\label{eq:TrueInnerProdEst}
\normtwo{\res_k}^2 - |\langle \res_k, \atom_{p_{k+1}}\rangle|^2 
  \leq \normtwo{\res_k}^2 - (1-\delta)^2|\cres_k(p_{k+1})|^2
  \leq \|\res_k\|_2^2.
 \end{equation}

%
 Using the triangle inequality and \eqref{eq:TrueInnerProdEst}, we obtain 
 \[
  \begin{split}
  \normtwo{\res_k - \cres_k(p_{k+1}) \atom_{p_{k+1}}}^2
  & \leq ( \normtwo{\res_k - \langle \res_k, \atom_{p_{k+1}}\rangle \atom_{p_{k+1}}} + \epsilon \sumcres_k)^2\\
  & =  \normtwo{\res_k - \langle \res_k, \atom_{p_{k+1}}\rangle \atom_{p_{k+1}}}^2
  + 2\epsilon \sumcres_k\normtwo{\res_k - \langle \res_k, \atom_{p_{k+1}}\rangle \atom_{p_{k+1}}} + \epsilon^2 \sumcres_k^2\\
  & \leq  \normtwo{\res_k}^2 \left(1 - \frac{(1-\delta)^2|\cres_k(p_{k+1})|^2}{\normtwo{\res_k}^2}+  \frac{\epsilon^2 \sumcres_k^2}{\normtwo{\res_k}^2} + \frac{2 \epsilon \sumcres_k}{\normtwo{\res_k}}\right).
  \end{split}
 \]

 Hence, 
 \[
  \frac{\epsilon^2 \sumcres_k^2 + 2 \epsilon \sumcres_k\normtwo{\res_k}}{(1-\delta)^2} 
  \leq 
  \frac{\delta^2 |\cres_k(p_{k+1})|^2 + 2 \delta |\cres_k(p_{k+1})|\normtwo{\res_k}}{(1-\delta)^2} 
  < |\cres_k(p_{k+1})|^2
 \]
ensures that $\|\res_{k+1}\|_2 < \|\res_k\|_2$. The latter inequality is equivalent to
\[
  \begin{split}
    \frac{2\delta}{1-2\delta}\normtwo{\res_k} & <  |\cres_k(p_{k+1})|.
  \end{split}
 \] \qed

\subsection{Proof of Theorem \ref{thm:Resdecrease2}}
%
It is clear that there exist $\delta$ and $\varepsilon$ such that $\frac{2\delta}{1-(2\delta +\varepsilon)} < \sqrt{\lambda_{\text{min}}}$. With such a choice, we have that 
\begin{equation}\label{eq:inprodest} \max_p |\langle \vect{y},\vect{d}_p\rangle |^2 \geq \lambda_{\text{min}}\| \vect{y}\|_2^2 > \left(\tfrac{2\delta}{1-(2\delta +\varepsilon)}\right)^2 \|\vect{y}\|_2^2.\end{equation} 

First, assume that $k=0$ or either of the stopping conditions in \eqref{eq:crit1} 
was met in the $(k+1)$-th selection step and $\cres_k$ is re-initialized with $\cres_k := \Dict^\ast (\sig - \Dict\coef_k) = \Dict^\ast \res^{\text{out}} = \Dict^\ast \res_k$. In other words, $\cres_k(p) = \langle \res_k,\atom_p\rangle$, for all $p$ and we have $k^{\text{out}} = k$. Furthermore, $p_{k+1} = \argmax_p |\cres_k(p)| = \argmax_p |\langle \res_k,\atom_p\rangle|$. With this, the right-hand side of the first inequality in \eqref{eq:crit1} equals $0$, such that it is trivially not satisfied. On the other hand, the second inequality in \eqref{eq:crit1} reduces to 
\[
 \max_p |\langle \res_k,\atom_p\rangle| 
 < \frac{2\delta (1-\varepsilon (1+\delta)^{-2}  \lambda_{\text{min}})^{\frac{k^{\text{out}}-k}{2}}}{1-(2\delta+\varepsilon)} \normtwo{\res^{\text{out}}}
 = \frac{2\delta}{1-(2\delta+\varepsilon)} \normtwo{\res_k} < \lambda_{\text{min}} \normtwo{\res_k},
\]
contradicting the definition of $\lambda_{\text{min}}$.

Hence, we can assume at every selection, that neither condition in \eqref{eq:crit1} 
is satisfied. Therefore, for every $k\geq 0$, and with $\sumcres_k:= \sum_{l=k^{\text{out}}+1}^{k} |\cres_{l-1}(p_{l})|$ (note that $k^{\text{out}}$ depends on $k$), we have 
\[
 (1+\delta)|\cres_k(p_{k+1})| 
 \geq |\cres_k(p_{k+1})| + \epsilon \sumcres_k 
 \geq \max_p |\langle \res_k, \atom_p\rangle| \geq \sqrt{\lambda_{\text{min}}}\|\res_k\|_2.
\]
Although $k^{\text{out}}$ depends on $k$, the estimate $(1+\delta)|\cres_k(p_{k+1})|\geq \sqrt{\lambda_{\text{min}}}\|\res_k\|_2$ does not. Clearly, for $l=0$, the inequality $|\cres_l(p_{l+1})| \geq \frac{2\delta}{1-(2\delta+\varepsilon)}\normtwo{\res_l}$ holds by \eqref{eq:inprodest}. Assume that it holds for for all $0\leq l < k$, for some fixed $k\geq 0$. Then, similar to the proof of Theorem \ref{thm:Resdecrease}
\[ 
  \begin{split}
  \epsilon^2 \sumcres_l^2 + 2\epsilon \sumcres_l\normtwo{\res_l} -  (1-\delta)^2|\cres_l(p_{l+1})|^2 
  & \leq \delta^2 |\cres_l(p_{l+1})|^2 + 2 \delta |\cres_l(p_{l+1})|\normtwo{\res_l} -  (1-\delta)^2|\cres_l(p_{l+1})|^2\\
  & \leq -\varepsilon |\cres_l(p_{l+1})|^2, 
  \end{split}
  \]
  such that 
 \[
    \normtwo{\res_{l+1}}^2 = \normtwo{\res_l - \cres_l(p_{l+1}) \atom_{p_{l+1}}}^2 
    \leq \left(1-\varepsilon\frac{|\cres_l(p_{l+1})|^2}{\normtwo{\res_l}^2}\right)\normtwo{\res_l}^2
    \leq \left(1-\varepsilon \frac{\lambda_{\text{min}}}{(1+\delta)^2}\right)\normtwo{\res_l}^2.
 \]
 With $\res_0 = \sig$, we obtain 
 \begin{equation}\label{eq:exponentialEst}
    \normtwo{\res_{k}} 
    \leq \left(1-\varepsilon  (1+\delta)^{-2} \lambda_{\text{min}}\right)^{\frac{k^{\text{out}}-k}{2}}\normtwo{\res^\text{{out}}} 
    \leq \left(1-\varepsilon  (1+\delta)^{-2} \lambda_{\text{min}}\right)^{k/2}\normtwo{\sig}.
 \end{equation}
 Since the second inequality in \eqref{eq:crit1} is not satisfied, the \eqref{eq:exponentialEst} implies 
 \[
 |\cres_k(p_{k+1})| \geq \tfrac{2\delta \normtwo{\res_k}}{1-(2\delta+\varepsilon)}, 
 \]
 completing the induction step and thus the proof.\qed

\subsection{Proof of Corollary \ref{cor:convergence}}
  Under the conditions of the corollary, we can invoke Theorem \ref{thm:Resdecrease2}. The result easily follows by considering \eqref{eq:exponentialEst} in the proof of Theorem \ref{thm:Resdecrease2} and noting that $k > \frac{\log(E)}{\log(1-\varepsilon (1+\delta)^{-2} \lambda_{\text{min}})}$ implies $(1-\varepsilon (1+\delta)^{-2} \lambda_{\text{min}})^k < E$. \qed

\end{document}